\documentclass[twocolumn]{article}
\usepackage{template/imr}
\usepackage{graphicx}
\usepackage{hyperref}
\usepackage{subcaption}
\usepackage{amssymb,bm}
\usepackage{algorithm2e}
\usepackage{graphicx}
\usepackage{amsmath}
\usepackage{xcolor}

\usepackage{tikz}

\begin{document}

\title{\uppercase{Naturally curved quadrilateral mesh generation using an adaptive spectral element solver}}

\author{Julian Marcon$^1$
  \and David A. Kopriva$^{2,3}$
  \and Spencer J. Sherwin$^4$
  \and Joaquim Peir\'{o}$^5$}
\date{
  $^1$Imperial College London, London, U.K. julian.marcon14@imperial.ac.uk\\
  $^2$The Florida State University, Tallahassee, FL, U.S.A. kopriva@math.fsu.edu\\
  $^3$San Diego State University, San Diego, CA, U.S.A.\\
  $^4$Imperial College London, London, U.K. s.sherwin@imperial.ac.uk\\
  $^5$Imperial College London, London, U.K. j.peiro@imperial.ac.uk\\
}

%% \abstract{*} and \keywords{*} must be before \maketitle.
\abstract{
% The number of words in the abstract should not exceed 150.

We describe an adaptive version of a method for generating valid naturally curved quadrilateral meshes.
The method uses a guiding field, derived from the concept of a cross field, to create block decompositions of multiply connected two dimensional domains.
The \textit{a priori} curved quadrilateral blocks can be further split into a finer high-order mesh as needed.
The guiding field is computed by a Laplace equation solver using a continuous Galerkin or discontinuous Galerkin spectral element formulation.
This operation is aided by using \textit{p}-adaptation to achieve {\color{black}faster convergence} of the solution with respect to the computational cost.
From the guiding field, irregular nodes and separatrices can be accurately located.
A first version of the code is implemented in the open source spectral element framework \emph{Nektar++} and its dedicated high order mesh generation platform \emph{NekMesh}.
}


\maketitle
\thispagestyle{empty}
\pagestyle{empty}

% Number of Pages should not exceed 18 pages including Figures and References.

\section{Introduction}

Finite element~\cite{Viertelabs-1708-02316}, finite volume, (block structured) finite difference~\cite{kozdon:2013} and spectral element methods~\cite{Karniadakis:2005fj,Kopriva:2009nx} all benefit from using quadrilateral meshes.
Yet fully automatic quadrilateral mesh generation is still difficult~\cite{Bommes:2013pd}, especially for high order (curved) elements used in spectral element methods.

A method was described in~\cite{Marcon} that generates large block structures with naturally curved boundaries that can be used as is, say for spectral element or block structured finite difference methods, or further subdivided into small quadrilaterals for finite element or finite volume methods.
The idea behind the approach came from the computer graphics community, which used cross fields~\cite{Viertelabs-1708-02316,Bommes:2013pd}.
The procedure in~\cite{Marcon} used a high order approach to compute a guiding field, accurately locating the irregular nodes, and accurately integrating the separatrices used to subdivide the domain.
It used a high order continuous or discontinuous Galerkin spectral element method (SEM or DGSEM) on a base triangular mesh, the approximation depending on the regularity of the boundary conditions needed to solve Laplace problems that determine the guiding field.
Unlike previous work, accuracy was preserved because all operations were performed on the guiding field rather than a cross field.
That new procedure was implemented in the open source spectral/\textit{hp} element framework \emph{Nektar++} and its high order mesh generation utility \emph{NekMesh}.

One of the significant costs in the high resolution field guided approach is the solution to the Laplace problems on the triangular mesh.
As the polynomial order of the approximation gets high, the condition number of the matrices that must be solved (either directly or iteratively) becomes large, leading to increased cost to get the field solution.
The approximation order was chosen \textit{ad hoc}, and there was no way to tell whether or not the chosen approximation order was sufficient.

One can avoid the \textit{ad hoc} nature of the choice of approximation order, and the slow convergence rates associated with choosing one that is higher than is necessary, by using an adaptive solver that adjusts the approximation to satisfy a desired error tolerance.
With a spectral element solver, adaptation can be done either by subdividing the underlying triangular elements (\textit{h}-adaptation) and/or by increasing the polynomial order (\textit{p}-adaptation).
In this work, we employ \textit{p}-adaptation because of its best accuracy, compared to \textit{h}-adaptation, with respect to the computational cost for smooth fields without discontinuities~\cite{Li2010}.
We use a dynamic \textit{p}-adaptive strategy previously devised for compressible flow simulation~\cite{Ekelschot2016,Moxey2017} and implemented in \emph{Nektar++}.

The methodology consists of three main steps.
First, a guiding field is computed.
This requires the generation of a high order triangular mesh of the domain and the solution to the Laplace problem using a \textit{p}-adaptive SEM or DGSEM solver.
Then, analysis of the guiding field must be carried out to construct a separatrix graph.
This construction process includes the location of interior critical points, the calculation of their valence as well as the valence of corners, and then the integration of streamlines and their possible manipulation.
Finally, a high order quadrilateral mesh can be generated.
To achieve this, streamlines are splined and used to split the CAD model into quadrilateral blocks, which can then be trivially meshed and refined.

This paper is structured as follows.
Section~\ref{sec:MathematicalFormulation} gives an outline of the mathematical formulation of the problem, including a definition of crosses, the guiding field and the partial differential equations to solve, a summary of the SEM/DGSEM formulation, and a overview of the \textit{p}-adaptation strategy.
Section~\ref{sec:implementation} covers some of the implementation details and Section~\ref{sec:examples} illustrates the methodology with examples.

\section{Formulation}\label{sec:MathematicalFormulation}

The goal is to subdivide any simply or multiply connected two dimensional domain, $\Omega$, with piecewise smooth boundaries into quadrilateral subdomains.
The decomposition is called regular if all the corners of the quadrilaterals have valence four, that is, each corner node connects four edges.
Especially in multiply connected domains, there will be \emph{irregular} nodes, where the valence will be larger or smaller than four.
The problem is how to generate a quadrilateral decomposition that has a minimum number of irregular nodes.

\subsection{The Guiding Field}\label{sec:guiding-field}

An approach that has been used in the computer graphics community to generate quadrilateral meshes over surfaces relies on the idea of a \emph{cross field}.
For full discussion of cross fields and methods based on them, see, for example, references~\cite{Viertelabs-1708-02316,Viertel:2017rt,Bunin2006,2008arXiv0802.2399B}.
We will present the basic idea.

A cross has four way rotational symmetry and can be represented by the four vectors
{\color{black}
\begin{equation}
C\left( \psi  \right) = \left\{ {{{\vec c}_k} = {{\left( { c_x, c_y} \right)}_k} = \left( {\cos \left( {\psi_k} \right),\sin \left( {\psi_k} \right)} \right)} \right\},
\label{eq:CrossDef}
\end{equation}
}
where
{\color{black}
\begin{equation}
\psi_k = \psi  + k \frac{\pi}{2}
\end{equation}
}
for, $k = 0,1,2,3$.
In turn, the tangent angle, or phase, $\psi$, can be computed from the four quadrant inverse tangent,
\begin{equation}
\psi = \frac{1}{4}\rm{atan2}\left({v,u}\right)\in\left[-\frac{\pi}{4},\frac{\pi}{4}\right].
\label{eq:PsiDef}
\end{equation}
In other words, a cross is represented by four unit vectors {\color{black}$\vec c_{k}$} at any point $\vec x \in \Omega$ in a vector field $\vec v = \left(u\left(\vec x\right),v\left(\vec x\right)\right)$ at which $\vec v\ne 0$.
A cross rotated by an angle $\pi/2$ is still the same cross.

The function $\psi$ is either tangent or orthogonal (due to the jump in the \(atan2\) function) to the streamlines of the field $\vec v$.
It is parallel to one of the branches of a cross and is undefined at a critical point of the field, i.e.\ where $\vec v = (0,0)$.
Given the way $\psi$ is computed from the arctangent, there will be a jump in value $\pi/2$ in $\psi$ depending on the signs of $u$ and $v$ even if $\vec v$ is smooth.
We will call lines in the field across which $\psi$ jumps \emph{jump lines}.

The method rests on finding an appropriate guiding field $\vec v$.
The only known values are on the boundary $\partial \Omega$ of $\Omega$.
There, the field is aligned with the boundary to ensure that the mesh to be generated is aligned there.
Along the boundary,
\begin{equation}
{{\vec v}_b} = \left(u_{b},v_{b}\right)=\left( {\cos \left(  4\theta_{b} \right),\sin ( 4\theta_{b})} \right),
\label{eq:BCs}
\end{equation}
where $\theta_{b}$ is the tangent angle of the boundary.
The vector $\vec v_{b}$ at any point along the boundary defines a cross~\eqref{eq:CrossDef} at that point.
The factor of four in eq.~\eqref{eq:BCs} ensures the same $u,v$ values for each $90^{\circ}$ rotation of the angle {\color{black}\(\theta_b\)}, and hence the rotational symmetry.

Boundary values defined with~\eqref{eq:BCs} are smoothly propagated to the interior.
Solving a Laplace problem for $u$ and $v$ 
\begin{equation}
\left\{ \begin{gathered}
  {\nabla ^2}\vec v = 0,\quad \, \vec x  \in \Omega  \hfill \\
  \vec v = {{\vec v}_b},\quad \, \vec x  \in \partial \Omega  \hfill \\ 
\end{gathered}  \right.
\label{eq:LaplaceEqns}
\end{equation}
guarantees that the field $\vec v$ is smooth in the interior of the domain.

\subsubsection{Critical Points}\label{sec:math_singularities}

To obtain a valid block decomposition of an arbitrary simply or multiply connected domain, \emph{irregular nodes} must be found where the valence \(\mathcal V \neq 4\).
In the cross field literature, these correspond to \textit{singular points} of the cross field where crosses are undetermined.
A similar phenomenon is observed in this guiding field method.
We return to Eq.~\eqref{eq:PsiDef} where \(\psi\) is undefined when \(u = v = 0\).
Points where \(\vec v = \vec 0\) are called \textit{critical points} of the guiding field and are a direct analogy of singular points of cross fields.

\begin{figure}[htb!]
  \centering
  \includegraphics[trim={15cm 0cm 15cm 10cm},clip,width=\linewidth]{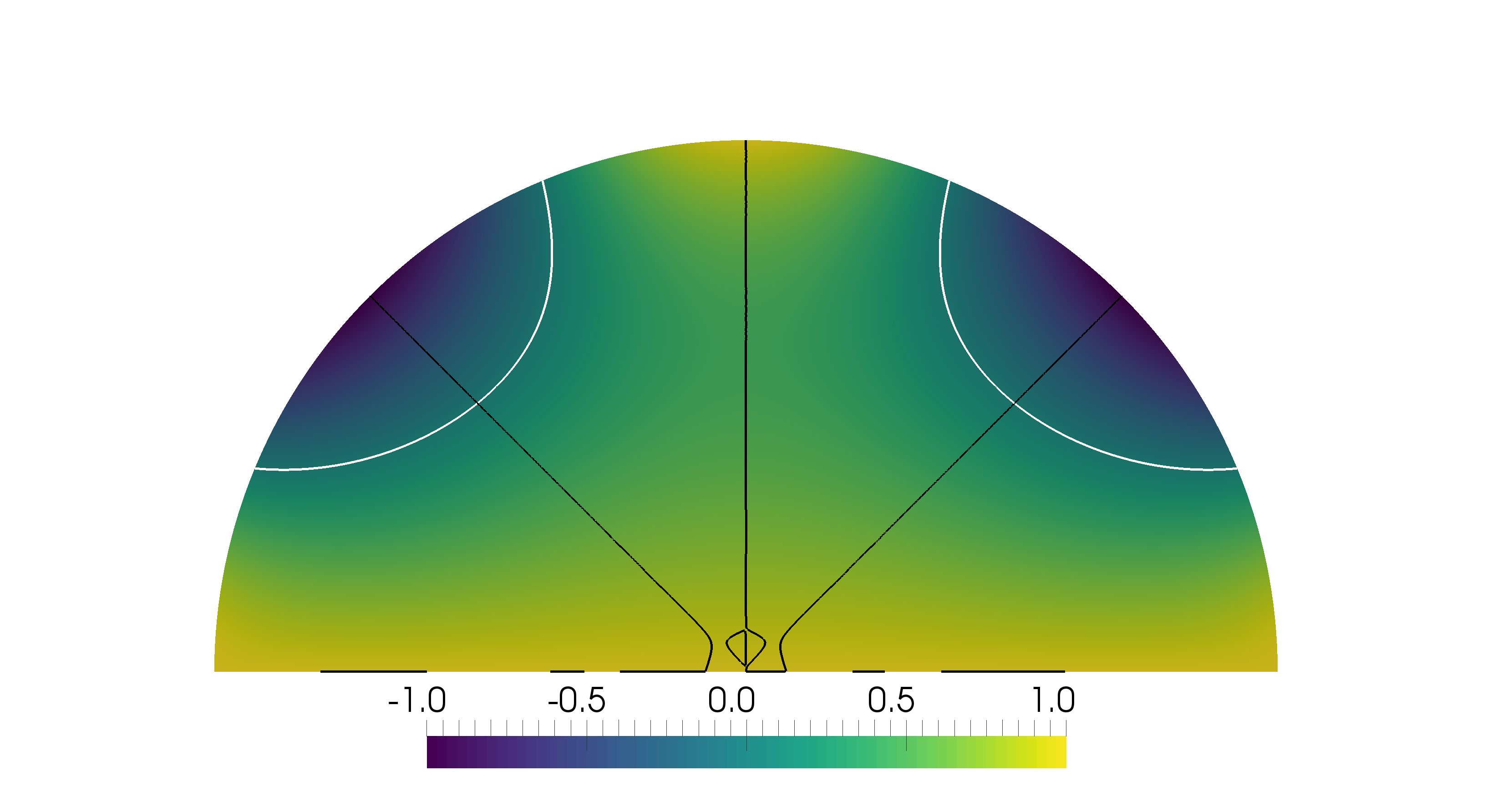}
  \includegraphics[trim={15cm 0cm 15cm 10cm},clip,width=\linewidth]{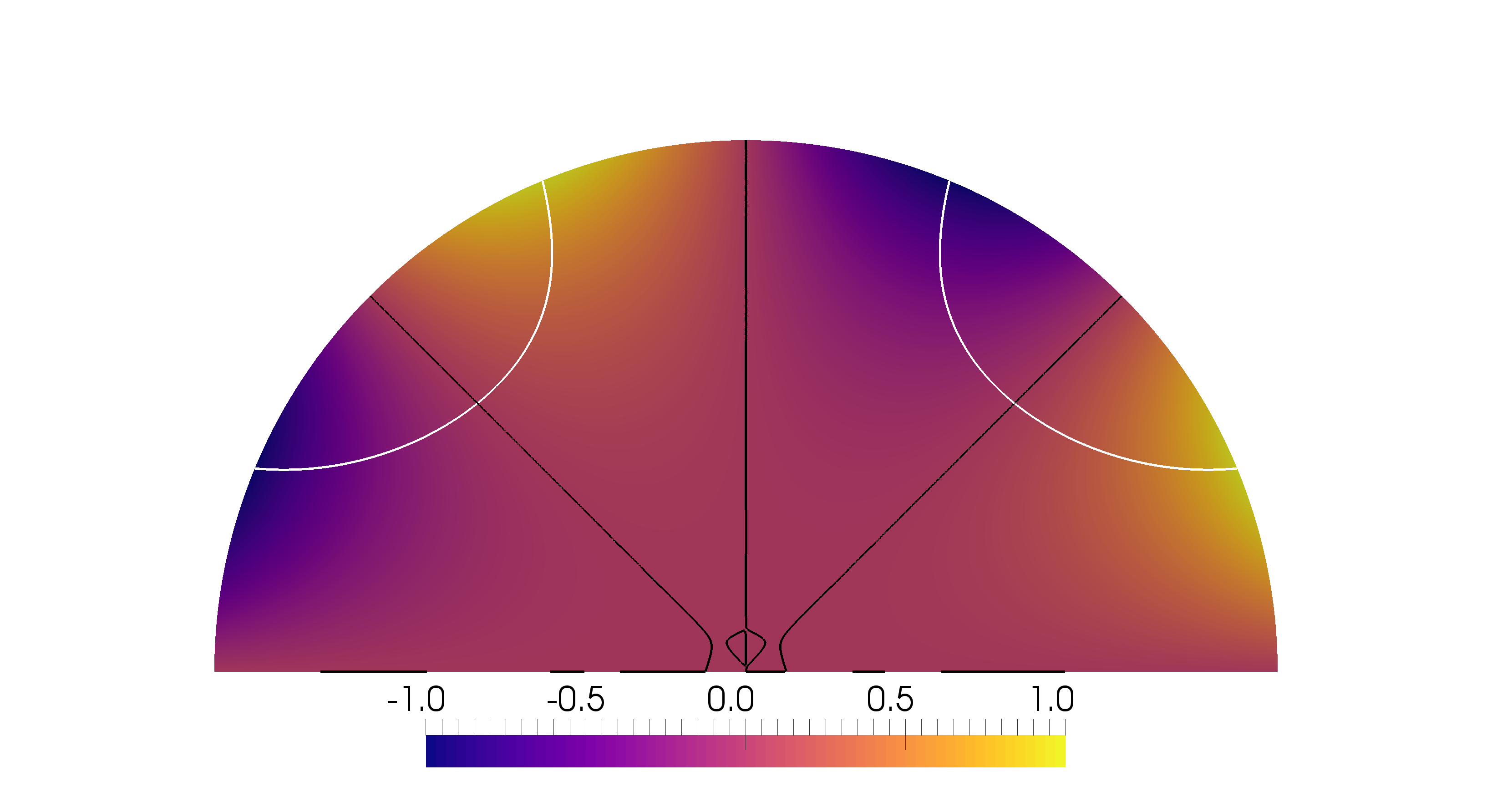}
  \includegraphics[trim={15cm 0cm 15cm 10cm},clip,width=\linewidth]{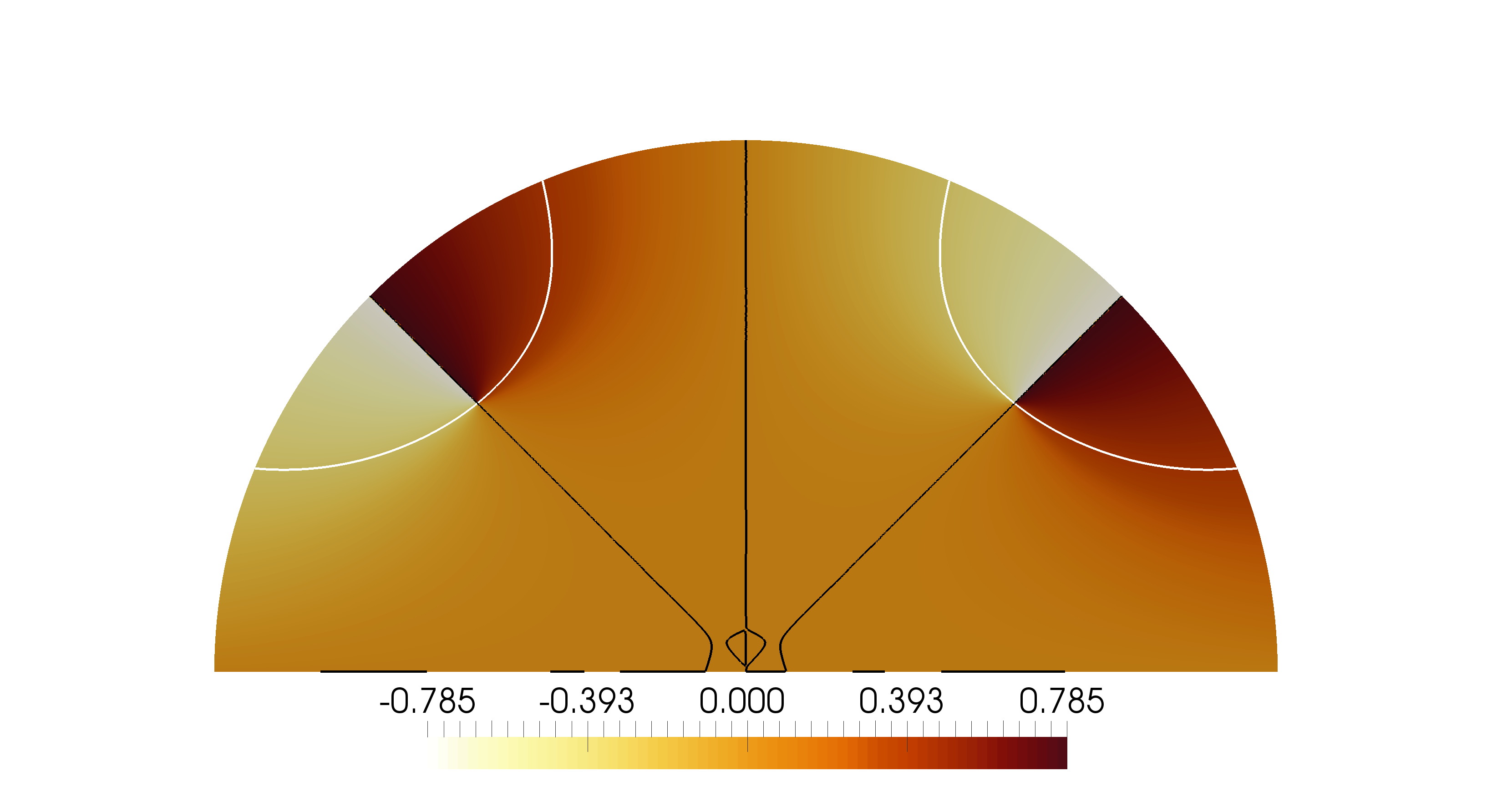}
  \caption{Solution \(u\) (top) and \(v\) (middle) fields and computed \(\psi \) field (bottom) on the half disc geometry. Isocontours of \(u = 0\) and \(v = 0\) shown in white and black respectively.}\label{fig:half-disc-solution}
\end{figure}

Critical points are tied to null isocontours in the domain.
Zeroes in \(u\) and \(v\) emerge from boundary conditions and propagate into the interior of the domain due to the Laplace equation.
This is shown in Fig.~\ref{fig:half-disc-solution} where null isocontours are shown intersecting at critical points.
Although \(\psi\) is undefined at critical points, \(\vec v\) remains smooth across the domain

We further need to calculate the valence of these irregular nodes, which will correspond to the number of quads connected to them.
We categorize them analytically, as recognizable topologies in the \(\vec v\) field~\cite{Tricoche:2001:CTS:601671.601695}, by the Poincar\'{e} index
\begin{equation}
  {i_\gamma } = \frac{1}{{2\pi }}\oint_\gamma  {d\phi },
\label{eq:pointcare}
\end{equation}
where {\color{black} \(\gamma\) is a simple arbitrary closed curve around a critical point and}
\begin{equation}
  \phi = \rm{atan2}\left({v,u}\right) = 4 \psi.
  \label{eq:phi-psi}
\end{equation}

A regular node in the domain would yield an index of zero while irregular nodes have an index of \(\pm 1\).
The index can, in turn, be used to calculate the number of separatrices connected to this irregular node, and therefore its valence.
The index \(i_{\gamma} = -1\) leads to \(\mathcal V = 5\) and \(i_{\gamma} = +1\) to \(\mathcal V = 3\)~\cite[Lemma 5.1]{Viertelabs-1708-02316},~\cite{Nicolas-Kowalski:2012fu}.

We adapt the Poincar\'{e} index~\eqref{eq:pointcare}, using relation~\eqref{eq:phi-psi} to the guiding field to obtain
\begin{equation}
  I_{c} = \frac{1}{\pi/2} \oint_c  {\frac{d\psi}{d\theta}d\theta },
  \label{eq:IDef}
\end{equation}
where $c$ is a (small) counter-clockwise circular contour centered on the critical point.
From this, the \emph{valence} \(\mathcal V\) of any point in the domain is defined as
\begin{equation}
  \mathcal V= 4 - I_{c}.
  \label{eq:singularity_valence}
\end{equation}
{\color{black}
The contour radius \(c\) must be chosen so that the closed curve encloses only one critical point.
If \(c\) encloses multiple critical points, an incorrect combined valence will be computed.
For example, two enclosed \(\mathcal V = 5\) critical points will produce a combined valence of 6, while two opposite (one \(\mathcal V = 5\) and one \(\mathcal V = 3\)) enclosed critical points will cancel out and produce a combined regular node (\(\mathcal V = 4\)).
In practice, critical points are far enough from each other that such problems do not arise, and one can always compute the radius around any critical point from the distance to the nearest of the remaining critical points.
}

Eq.~\eqref{eq:singularity_valence} can be generalized to an open line circle integral to handle boundary points too~~\cite{2008arXiv0802.2399B}, leading to
\begin{equation}
  {I\left(\theta_{0},\theta_{f}\right) } = \frac{1}{\pi/2} \int_{{\theta _0}}^{{\theta _f}} {\frac{{d\psi }}{{d\theta }}d\theta },
\end{equation}
such that \(I_{c}=I\left(0,2\pi\right)\).
Similarly, the valence of any boundary point is defined as
\begin{equation}
  \mathcal V = {\frac{\Delta\theta}{\pi/2} - {I\left(\theta_{0},\theta_{f}\right) } } ,
  \label{eq:math_corner}
\end{equation}
where \( \Delta\theta = \theta_f - \theta_0 \).
\( I\left(\theta_{0},\theta_{f}\right) \) can be viewed as a correcting factor for the valence towards an integer value.

\subsubsection{Streamline Integration}

\emph{Separatrices} of the block decomposition, in the form of \emph{streamlines} of the guiding field, must be traced from irregular nodes and boundary corners throughout the domain.
The guiding field, akin to cross fields, can be used for this purpose too.
We define a streamline as \( \vec x (t)\) such as
\begin{equation}
  \frac{d \vec x }{dt} = \tilde v\left(\psi'( \vec x )\right)
  \label{eq:streamline-integration}
\end{equation}
where \(t\) is the integration parameter and \(\tilde v\left(\psi'( \vec x )\right)\) is the adjusted guiding field vector.

Because \(\psi\) is discontinuous, defined on \(\left[-\frac{\pi}{4},\frac{\pi}{4}\right]\) by ~\eqref{eq:PsiDef}, whereas crosses are invariant of rotations of \(\pi/2\), we use an \emph{adjusted} guiding direction \(\psi'\).
The guiding direction must indeed be adjusted by a multiple of \(\pi/2\) every time a streamline crosses a jump line, defined in Sec.~\ref{sec:guiding-field}.
This allows us to retain the consistent directionality of a smooth cross field.

Eq.~\eqref{eq:streamline-integration} can be integrated with respect to \(t\) given appropriate initial conditions
\begin{equation}
  \vec x (t_0) =  \vec x _0,
\end{equation}
where \(t_0\) is the initial integration time and \( \vec x _0\) the location of the start of the streamline.
The start of the streamline may either be an irregular node or a boundary corner.
The ordinary differential equation~\eqref{eq:streamline-integration} can be integrated using any traditional numerical scheme.
Sec.~\ref{sec:streamlines} will more specifically cover the integration strategy.

\subsection{The Spectral Element Method}

We compute the BVPs~\eqref{eq:LaplaceEqns} with either a continuous Galerkin (CG) or a discontinuous Galerkin (DG) spectral element method on a triangular mesh~\cite{Karniadakis:2005fj}.
Spectral element methods are spectrally accurate, meaning that the convergence rate depends only on the smoothness of the solution.
They are high resolution in that they use a large number of degrees of freedom within an element.
Unlike traditional low order finite element methods, the high order polynomial expansion of the solution inside each element allows us to locate critical points, identify their valence, and finally trace streamlines with high order accuracy.

In this work, we use the spectral/\textit{hp} element methods formulation described in detail in reference~\cite{Karniadakis:2005fj} and implemented in \emph{Nektar++}~\cite{Cantwell2015}.
We briefly survey the fundamentals of this method here for the numerical solution of PDEs, \(\mathcal{L} u = 0\), over a domain \( \Omega \).
The domain \(\Omega \) is covered by a mesh of triangular finite elements, \(\Omega_e\) so that \(\Omega = \cup \Omega_e\) and \(\Omega_{e_1} \cap \Omega_{e_2} = \partial\Omega_{e_1 e_2}\) is either an empty set or the interface between two elements, and is of one dimension less than the mesh.
We solve the PDE problem weakly and require that \(\left.u\right|_{\Omega_e}\) is in the Sobolev space {\color{black}\(H^{k}(\Omega_e)\) $k=0,1$}.
In the CG formulation, we require the solution to be in $H^{1}$, $H^{0}$ for DG.

We formulate the solution to~\eqref{eq:LaplaceEqns} in weak form: \(find\ \vec v \in H^k (\Omega)\ such\ that\)
\begin{equation}
  a(\vec v, \vec w) = l(\vec w) \quad \forall \vec w \in H^k (\Omega),
  \label{eq:WeakForm}
\end{equation}
where \(a(\cdot,\cdot)\) is a symmetric bilinear form, \(l(\cdot)\) is a linear form and \(H^k (\Omega)\) is the 
Sobolev space formally defined as
{\color{black}
\begin{equation}
  H^k (\Omega) := \{ \vec w \in L^2 (\Omega) \mid D^\alpha \vec w \in L^2 (\Omega)\ \forall\ \left|\alpha\right| \leq k \}.
\end{equation}
}

We solve problem \eqref{eq:WeakForm} numerically and therefore consider solutions in a finite dimensional subspace \(V_N \subset H^k(\Omega)\). In the finite dimensional subspace, the problem becomes: \(find\ \vec v^\delta \in V_N\ such\ that\)
\begin{equation}
  a(\vec v^\delta, \vec w^\delta) = l(\vec w^\delta) \quad \forall \vec w^\delta \in V_N,
\end{equation}
augmented with boundary conditions,~\eqref{eq:BCs}.
In the CG formulation, we also enforce the condition \(V_N \subset C^0\).

To represent \(\vec v^\delta( \vec x ) = \sum_n \hat{u}_n \Phi_n( \vec x )\), we take a weighted sum of \(N\) trial functions \(\Phi_n( \vec x )\) defined on \(\Omega \).
This transforms the problem of finding the subspace solution to one of finding the coefficients \(\hat{u}_n\) that define \(\vec v^\delta( \vec x )\) within an element.
To obtain a unique choice of coefficients \(\hat{u}_n\), we require that the residual \(R = \mathcal{L} \vec v^\delta \) be orthogonal to the polynomial space, i.e. that its \(L^2\) inner product, with respect to the test functions \(\Psi_n( \vec x )\), is zero.
In the Galerkin approximation, the projection the test functions are chosen to be the same as the trial functions, i.e. \(\Psi_n = \Phi_n\).

The contributions of each element in the domain must be combined to construct the global basis \(\Phi_n\).
A parametric mapping {\color{black}\(\mathcal{X}_e : \mathcal{E}\to \Omega_e\) exists to each element \(\Omega_e\) from a standard reference element \(\mathcal{E} \subset {[-1,1]}^d\)}.
This mapping is given by \( \vec x  = \mathcal{X}_e\left( \vec \xi\right)\),
where \( \vec x \) is the vector of the physical coordinates and \( \vec \xi \) defines the coordinates in the reference space.
In \emph{Nektar++}, triangular, tetrahedral, prismatic, and pyramidal elements are created by collapsing one or more of the coordinate directions to create singular vertices.
This allows it to support, for this work, more easily generated triangular meshes.

We use a local polynomial basis on the reference element to represent the solution.
A one-dimensional order-\(P\) basis is a set of polynomials \(\Phi_p(\xi), 0 \leq p \leq P\), defined on the reference segment \(-1 \le \xi_1 \le 1\).
A tensor basis is used in two dimensional reference regions, where the polynomial space is constructed as the tensor product of one dimensional bases on segments and quadrilaterals.

Finally, the discrete solution in a physical element \(\Omega_e\) can be expressed as
\begin{equation}
\vec v  ^\delta( \vec x ) = \sum_n \hat{v}_n \phi_n (\mathcal{X}_e^{-1} ( \vec x )),
\label{eq:SolutionInterpolant}
\end{equation}
with \(\hat{v}_n\) being the coefficients computed by the Galerkin procedure.
We therefore restrict the solution space to
\begin{equation}
  V := \{\vec v \in H^k(\Omega) \mid \left.\vec v\right|_{\Omega_e} \in \mathcal{P}_P(\Omega_e)\},
\end{equation}
where \(\mathcal{P}_P(\Omega_e)\) is the space of order \(P\) polynomials on \(\Omega_e\).

An assembly operator assembles the element contributions to the global solution.
In the CG formulation, elemental contributions of neighbours are summed to enforce \(C^0\)-continuity.
In the DG formulation, flux values are transferred from the element interfaces into the global solution vector.
{\color{black}
\emph{Nektar++} supports both discretizations, allowing us to trivially switch from one to the other at any time.
The DG formulation is necessary for the method to work with arbitrary corner angles in the geometry, as will be explained in Sec.~\ref{sec:solution}.
}

\subsubsection{\textit{p}-Adaptation}\label{sec:p-adapt}

Spectral element methods also offer the freedom to locally increase or decrease the resolution of the solution according to an estimated local error.
Adaptation is achieved by adding or removing modes --- i.e.\ increasing or decreasing the polynomial order of individual elements.
The \textit{p}-adaptive strategy employed in this paper follows the procedure laid out in references~\cite{Ekelschot2016} and~\cite{Moxey2017}.
The idea is simple: where local error is high, we increase local resolution; where local error is low, we decrease resolution.
The procedure is described in Algorithm~\ref{alg:p-adapt}, where \(e\) is an individual element, \(S_e\) and \(P_e\) are its associated error indicator and polynomial order, \(\epsilon_u\) and \(\epsilon_l\) are the upper and lower error thresholds, and {\color{black}\(P_{max}\)} and {\color{black}\(P_{min}\)} are the maximum and minimum polynomial orders allowed.

\begin{algorithm}[htb]
  \caption{The \textit{p}-adaptive procedure.}\label{alg:p-adapt}

  \While{any \(P_e\) is modified}{
    Calculate solution to~\eqref{eq:LaplaceEqns} \;
    \ForEach{\(e\)}{
      Calculate \(S_e\) \;
      \uIf{\(S_e>\epsilon_u \wedge P_e<{\color{black}P_{max}}\)}{
        increment \(P_e\) \;
      }
      \uElseIf{\(S_e<\epsilon_l \wedge P_e>{\color{black}P_{min}}\)}{
        decrement \(P_e\) \;
      }
      \Else{
        maintain \(P_e\) \;
      }
    }
  }
\end{algorithm}

We use an estimate of the discretization error from the formulation described in reference~\cite{Persson2006} for the error indicator (sometimes called \emph{sensor}).
{\color{black}
We only require to compute the solution \(u_P\) once at order \(P\), which we then project onto polynomial space \(P-1\) to obtain \(u_{P-1}\).
From these solutions at different polynomial orders, we calculate our sensor:
}
\begin{equation}
  S_e = \frac{ \left\Vert u_P - u_{P-1} \right\Vert^{2}_{2,e} }{ \left\Vert u_P \right\Vert^{2}_{2,e} },
\end{equation}
where \(\left\Vert \cdot \right\Vert_2\) is the \(L_2\)-norm.
For a modal expansion, this sensor is cheap to compute.
It phyically represents the ratio between the energy of the highest mode and the rest of modal energy spectrum.
The sensor \(S_e\) can be computed for either component of the solution \(\vec v\) in~\eqref{eq:LaplaceEqns}.

\section{Implementation of the Procedure}\label{sec:implementation}

The process is split into three main steps, each associated with different libraries and utilities of the \emph{Nektar++} platform~\cite{Cantwell2015,Moxey2019}:
\begin{enumerate}
  \item Computation of a guiding field by solving the Laplace problem~\eqref{eq:LaplaceEqns}. Carried out in the \emph{solver}.
  \item Analysis of the guiding field, including locating critical points, calculating their valence and the valence of corners, and the integration of streamlines. Carried out in the post-processing utility \emph{FieldConvert}.
  \item Generation of a quadrilateral mesh by splitting the CAD model into blocks and splitting these further into a fine mesh. Carried out in the meshing utility \emph{NekMesh}.
\end{enumerate}

\begin{figure}[htb!]
  \centering
  \includegraphics[trim={5cm 0 5cm 0},clip,width=\linewidth]{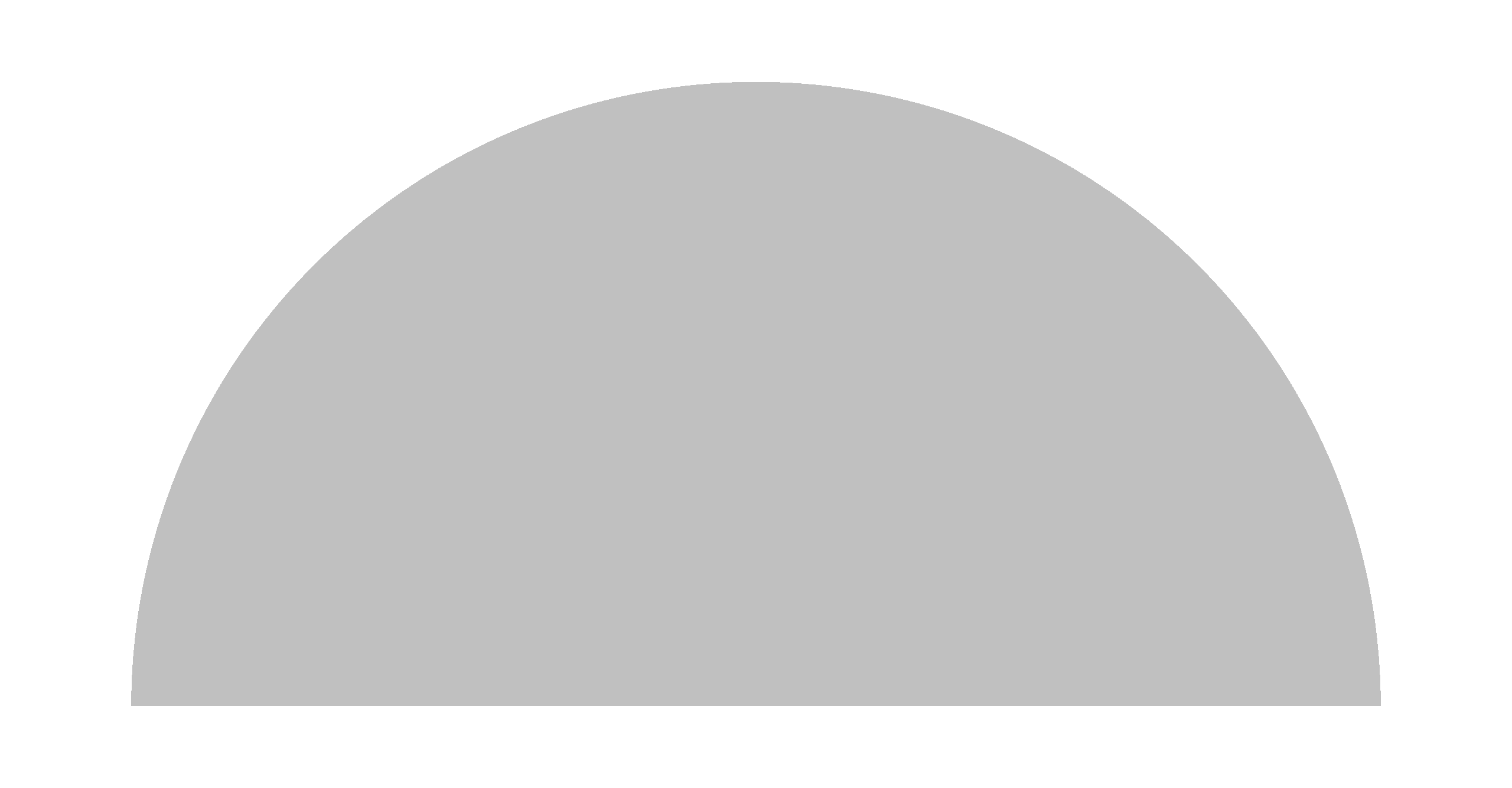}
  \caption{Geometry of a half disc for illustration of the quadrilateral meshing method.}\label{fig:half-disc-geometry}
\end{figure}

This section provides an overview of the process, illustrated on the geometry of a half disc geometry (see Fig.~\ref{fig:half-disc-geometry}), commonly used in the cross field literature~\cite{Nicolas-Kowalski:2012fu,Viertel:2017rt}.
The reader should refer to reference~\cite{Marcon} for implementation details.

\subsection{Computation of a Guiding Field}\label{sec:solution}

\begin{figure}[htb!]
  \centering
  \includegraphics[trim={5cm 0 5cm 0},clip,width=\linewidth]{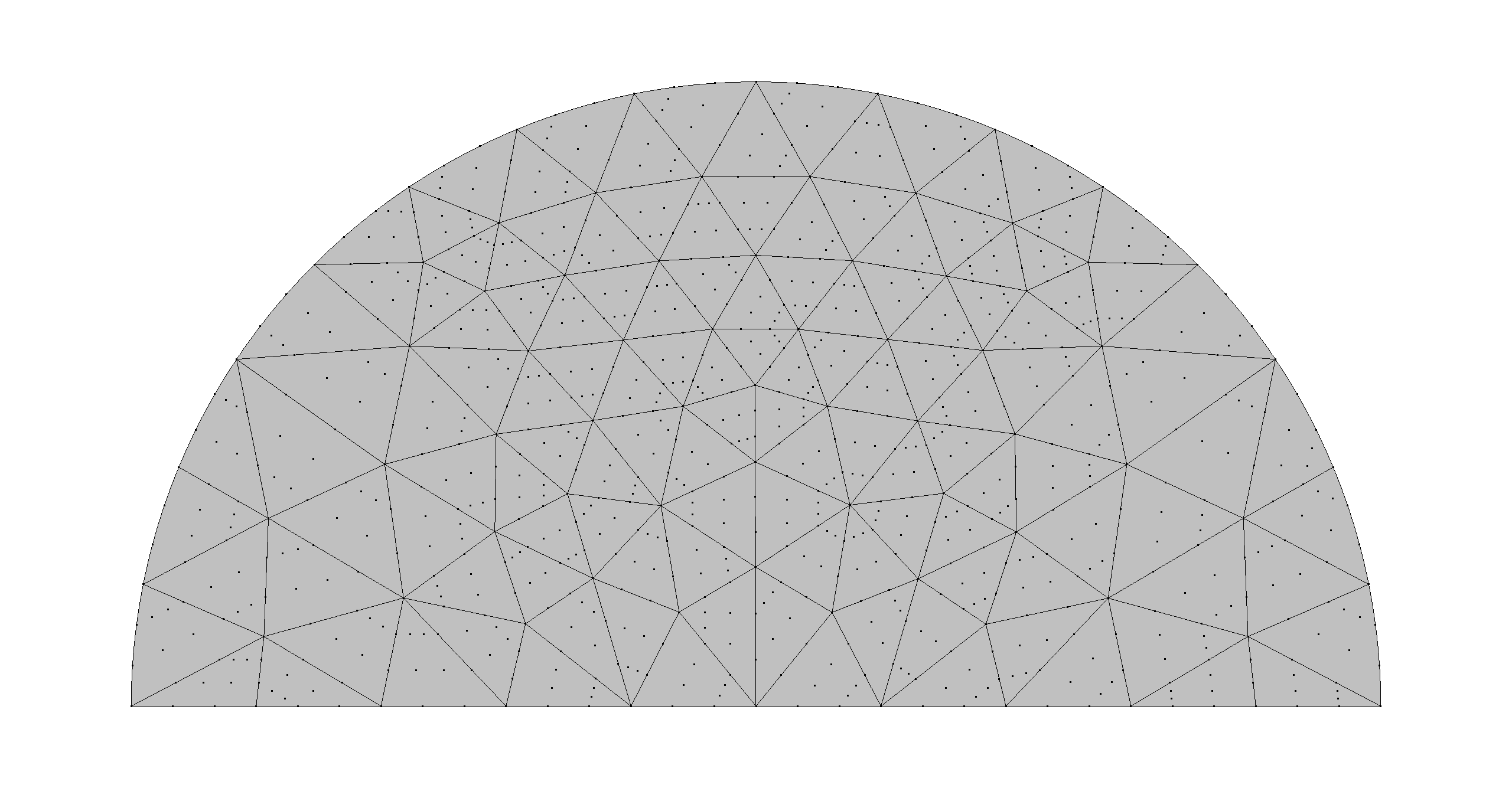}
  \caption{Third order triangular spectral element mesh for the half disc geometry.}\label{fig:half-disc-tris}
\end{figure}

A standard Laplace solver for the boundary value problem~\eqref{eq:LaplaceEqns} is used in the open source spectral/\textit{hp} element framework \emph{Nektar++}~\cite{Cantwell2015,Moxey2019} with boundary conditions~\eqref{eq:BCs}.
The geometry is passed onto the guiding field solver in the form of a high order triangular mesh previously generated in \emph{NekMesh}~\cite{Moxey2019} following the procedure laid out in reference~\cite{Sherwin2002}.
For the reference geometry, such a mesh is illustrated in Fig.~\ref{fig:half-disc-tris}.

For performance purposes, the triangular mesh should be coarse for fast computation of the guiding field.
Resolution will be added in areas of interest by use of \textit{p}-adaptation, as described in Sec.~\ref{sec:p-adapt}.
A coarse mesh also proved to be lighter to process in the analysis of the guiding field, as will seen in Sec.~\ref{sec:analysis}.

Notice that the boundary conditions (BCs)~\eqref{eq:BCs} depend on the geometry and may therefore be discontinuous where the geometry is not smooth.
To be exact, BCs are continuous if the boundaries satisfy one of two conditions:
\begin{itemize}
  \item The curve at each point is {\color{black}\(C^1\)-continuous}, i.e.~a smooth curve;
    or 
  \item The boundary curve is only {\color{black}\(C^0\)-continuous}, i.e.~a corner, and the angle is a multiple of \(\pi/2\).
\end{itemize}
Otherwise, they are discontinuous.

A discontinuous Galerkin formulation, available in \emph{Nektar++}, offers a discretization consistent way of handling discontinuous BCs without the need for \textit{ad hoc} smoothing of a corner's BCs as used in traditional low order techniques~\cite{Viertelabs-1708-02316}.

\begin{figure}[htb!]
  \centering
  \includegraphics[trim={15cm 0cm 15cm 5cm},clip,width=\linewidth]{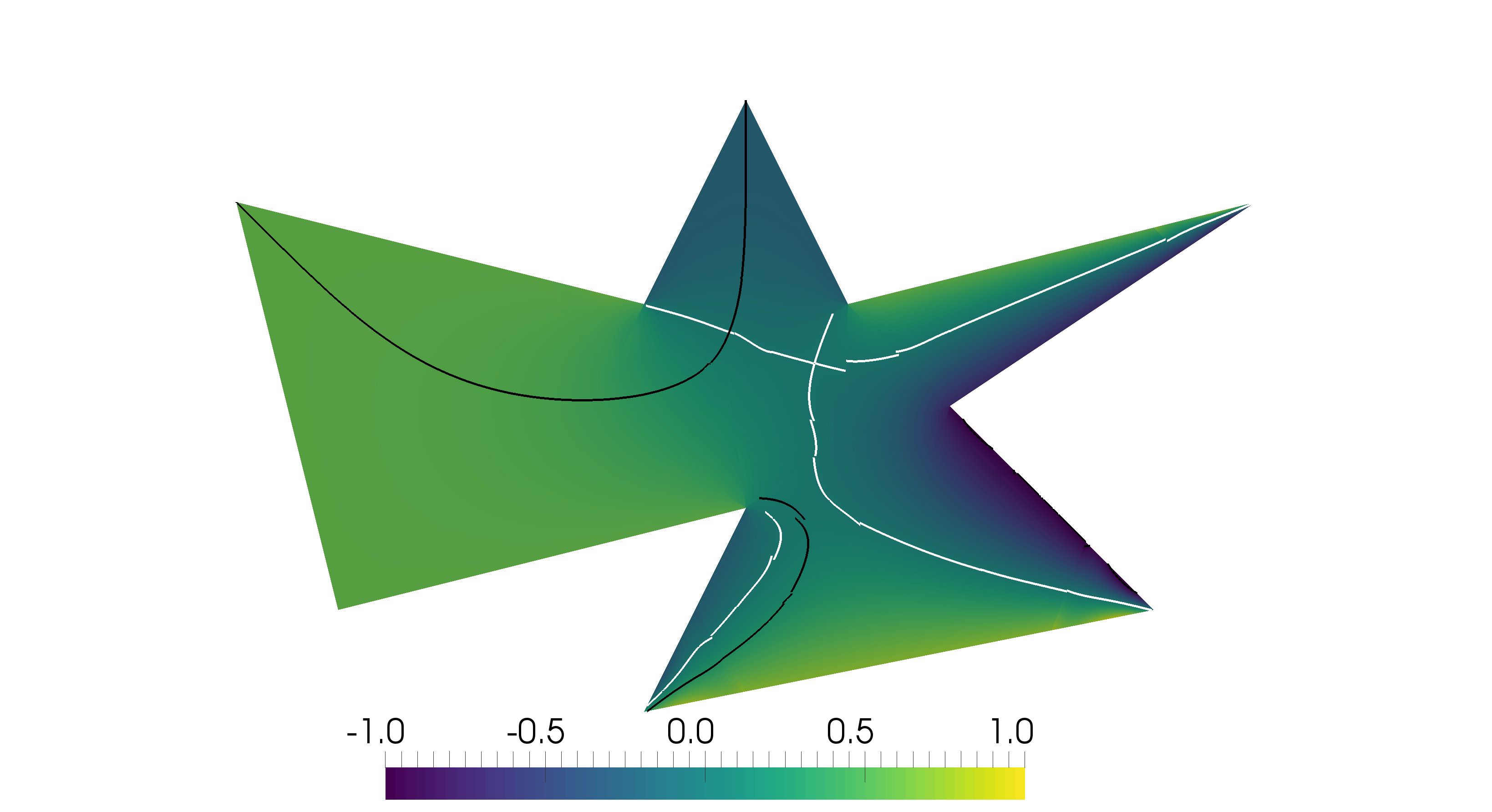}
  \includegraphics[trim={15cm 0cm 15cm 5cm},clip,width=\linewidth]{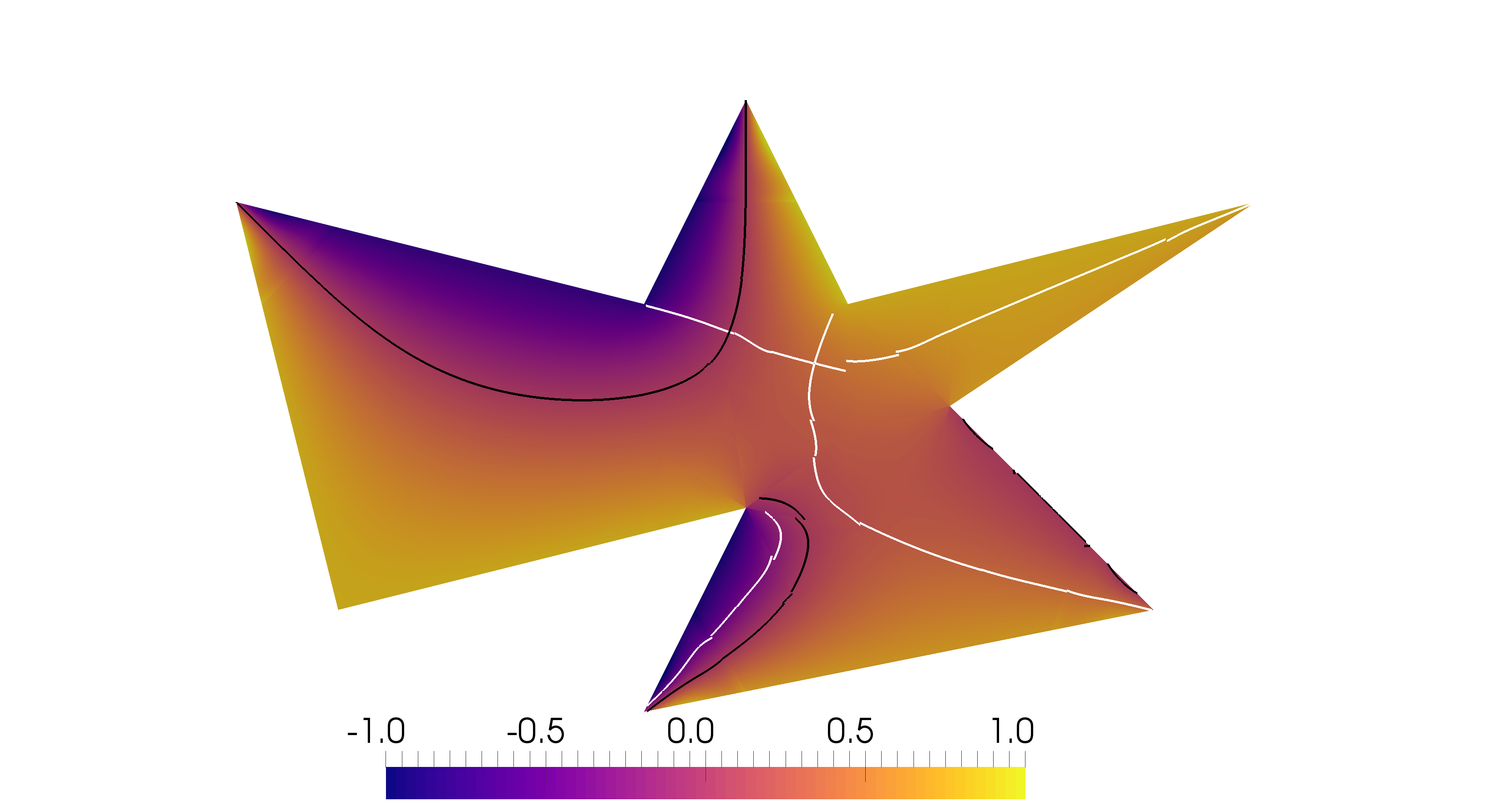}
  \includegraphics[trim={15cm 0cm 15cm 5cm},clip,width=\linewidth]{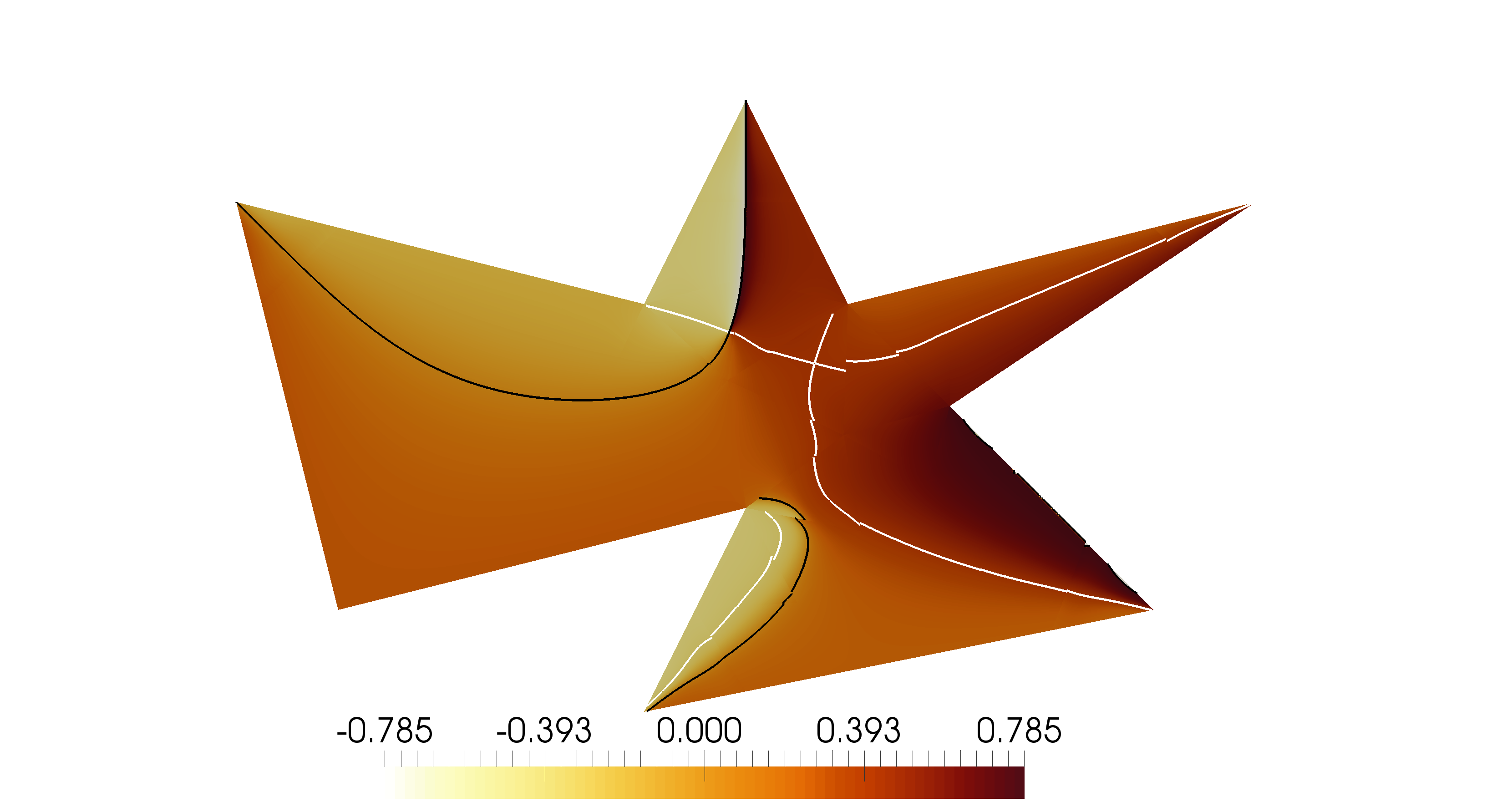}
  \caption{Solution \(u\) (top) and \(v\) (middle) fields and computed \(\psi \) field (bottom) on a polygon with acute and obtuse corners using a DG discretization. Isocontours of \(u = 0\) and \(v = 0\) are shown in white and black respectively.}\label{fig:polygon-solution}
\end{figure}

As an example of a geometry that requires a DG solve, we present the polygon geometry shown in Fig.~\ref{fig:polygon-solution}, which shows the solution and $\psi$ fields.
We note that the discontinuous BCs are naturally imposed with the Dirichlet BCs enforced through fluxes when needed.
Considerations beyond the scope of this paper are made to address the discontinuous nature of the field solution, but they however do not affect the rest of the procedure.

\begin{figure}[htb!]
  \centering
  \includegraphics[trim={15cm 2cm 15cm 10cm},clip,width=\linewidth]{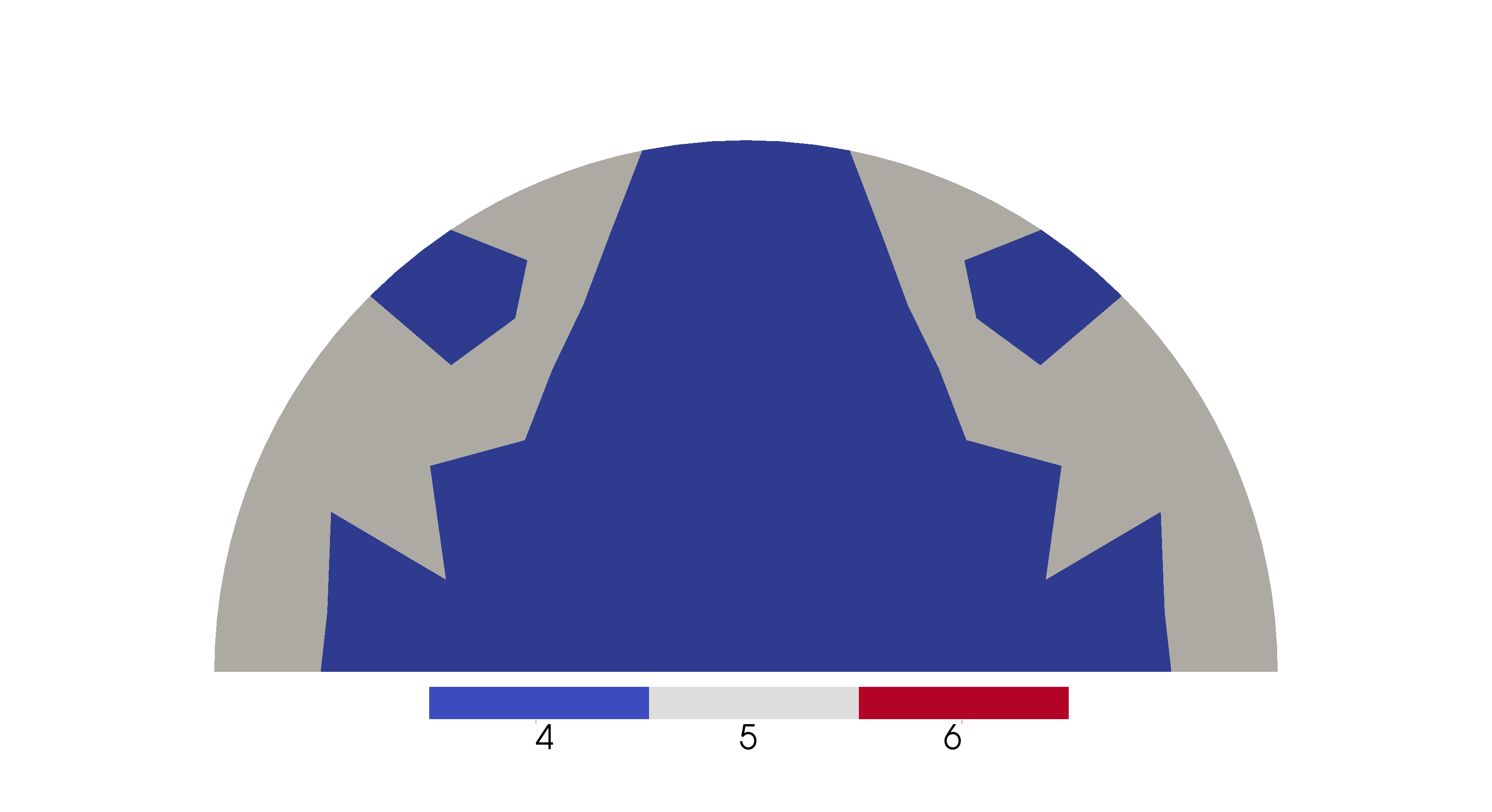}
  \includegraphics[trim={15cm 2cm 15cm 10cm},clip,width=\linewidth]{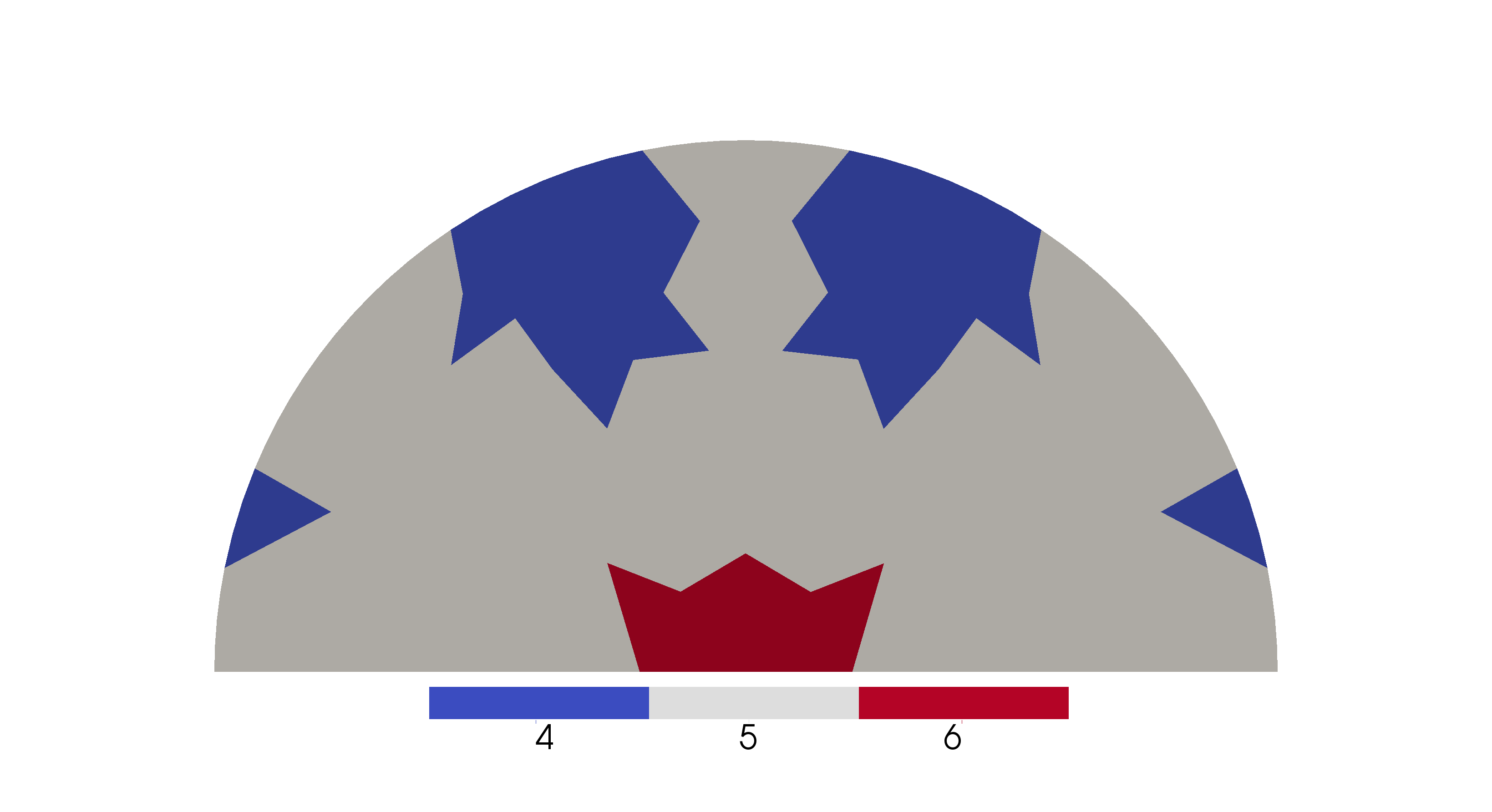}
  \caption{Distribution of local number of modes (\(=P+1\)) for half disc geometry. Adaptation based on \(u\) (top) and \(v\) (bottom).}\label{fig:half-disc-p}
\end{figure}

Following the calculation of a solution, \textit{p}-adaptation can be applied iteratively to locally increase or decrease the local polynomial order inside each element.
Using the approach described in Sec.~\ref{sec:p-adapt}, a finer solution is obtained without the need to uniformly increase the resolution and therefore the cost of the calculation.
Fig.~\ref{fig:half-disc-p} shows the final distribution of local number of modes when adapting based on either \(u\) (top) or \(v\) (bottom).
Higher polonomial orders are found near the boundaries where the gradient of the solution is stronger.

\subsection{Analysis of the Guiding Field}\label{sec:analysis}

The mesh generation procedure continues with an analysis of the guiding field, in a way where crosses are never generated.
The analysis is carried out in \emph{FieldConvert}, \emph{Nektar++}'s post-processing utility~\cite{Moxey2019}, and consists of three main parts:
\begin{enumerate}
  \item Location of interior critical points.
  \item Calculation of the valence of interior critical points and corners.
  \item Integration of streamlines.
\end{enumerate}
The goal of this analysis is to construct a separatrix graph to pass on to \emph{NekMesh} for the generation of a coarse high order quadrilateral mesh.

\subsubsection{Location of Interior Critical Points}

As described in Sec.~\ref{sec:math_singularities}, interior critical points are located at \(\vec v = \vec 0\).
To locate all of them we take a simple narrowing down approach.
We first identify elements that are likely to contain a \(\vec v = \vec 0\) point.
These elements should have at least one quadrature point with each of the values \( u \leq 0 \), \( u \geq 0 \), \( v \leq 0 \) and \( v \geq 0 \).
These elements either contain a critical point or are close to an element that does.
For each of these elements, we perform a Newton search for a point satisfying \(\vec v = \vec 0\).
A critical point may be found outside of the element boundaries, in which case we ignore it.
Fig.~\ref{fig:half-disc-solution} shows the location of the two interior critical points of the reference geometry as the crossing of the $u$ and $v$ zero contours, consistent with the literature~\cite{Nicolas-Kowalski:2012fu}.

We accelerate the search inside each element by performing Netwon's method in reference space.
Along with the initial identification of likely elements, all operations are performed in reference space at negligible computational cost.

\subsubsection{Calculation of Valences}

After the interior critical points are located, valences of these as well as of corners can be computed using~\eqref{eq:singularity_valence} and~\eqref{eq:math_corner} respectively.

Here again, crosses are not generated and valences can be computed from the guiding field alone.
We note that the closed curve circle integral \(I_c\) in~\eqref{eq:IDef} evaluates to zero when \(\psi\) is continuous.
This is the case of regular nodes in the domain.
Where a critical point appears, lines of \(\psi\) jumps appear and \(I_c\) can be evaluated based on the number and direction of these.
In a counter-clockwise manner, a positive jump (i.e.~from \(-\frac{\pi}{4}\) to \(\frac{\pi}{4}\)) leads to a negative {\color{black}index \( I_c = -1 \)}, and vice versa.
The sign of the jump can itself be determined solely by values of \(\vec v\), without the need to compute \(\psi \) or the cross field.

Using this knowledge, we can calculate the valence of interior critical points from~\eqref{eq:singularity_valence} by marching on a circle of small radius \(c\) around each critical point.
The same reasoning is applied to the calculation of corner valences with~\eqref{eq:math_corner} where \(\theta_0\) and \(\theta_f\) must also be taken into account and are found trivially.
{\color{black}
From this logic, we find that the two critical points of the reference geometry are characterized by a single negative jump line each, and therefore have a valence \(\mathcal V = 3\).
Likewise, the two corners have a trivially evaluated valence \(\mathcal V = 1\) due to \(\Delta\theta = \pi/2\) and the absence of any jump line in \(\psi\).
See Fig.\ref{fig:half-disc-solution}~(bottom) for reference.
}

Although these evaluations are performed in physical space at a negligible computational cost, they can further be accelerated by performing them in reference space.

\subsubsection{Integration of Streamlines}\label{sec:streamlines}

Streamlines emanate from the interior irregular nodes and corners.
Because \( \vec v = \vec 0 \) at irregular nodes, \(\psi \) is undefined.
Likewise, \(\psi\) at corners may be ambiguous if BCs are discontinuous.
It is therefore necessary to identify the initial direction of each streamline by looking at the value of \(\psi\) nearby the irregular nodes, at a small distance \(c\) from them.
The initial direction can be found iteratively by progressively refining an initial guess, akin to a fixed point search.
For the first streamline of an interior irregular node, any initial guess may be taken.
For subsequent streamlines as well as for streamlines coming from corners, an initial guess can be take at angles multiples of \( \frac{\Delta\theta}{\mathcal V} \) where \(\Delta\theta = 2\pi\) for an interior irregular node.

\begin{figure}[htb!]
  \centering
  \includegraphics[trim={30cm 0 30cm 0},clip,width=\linewidth]{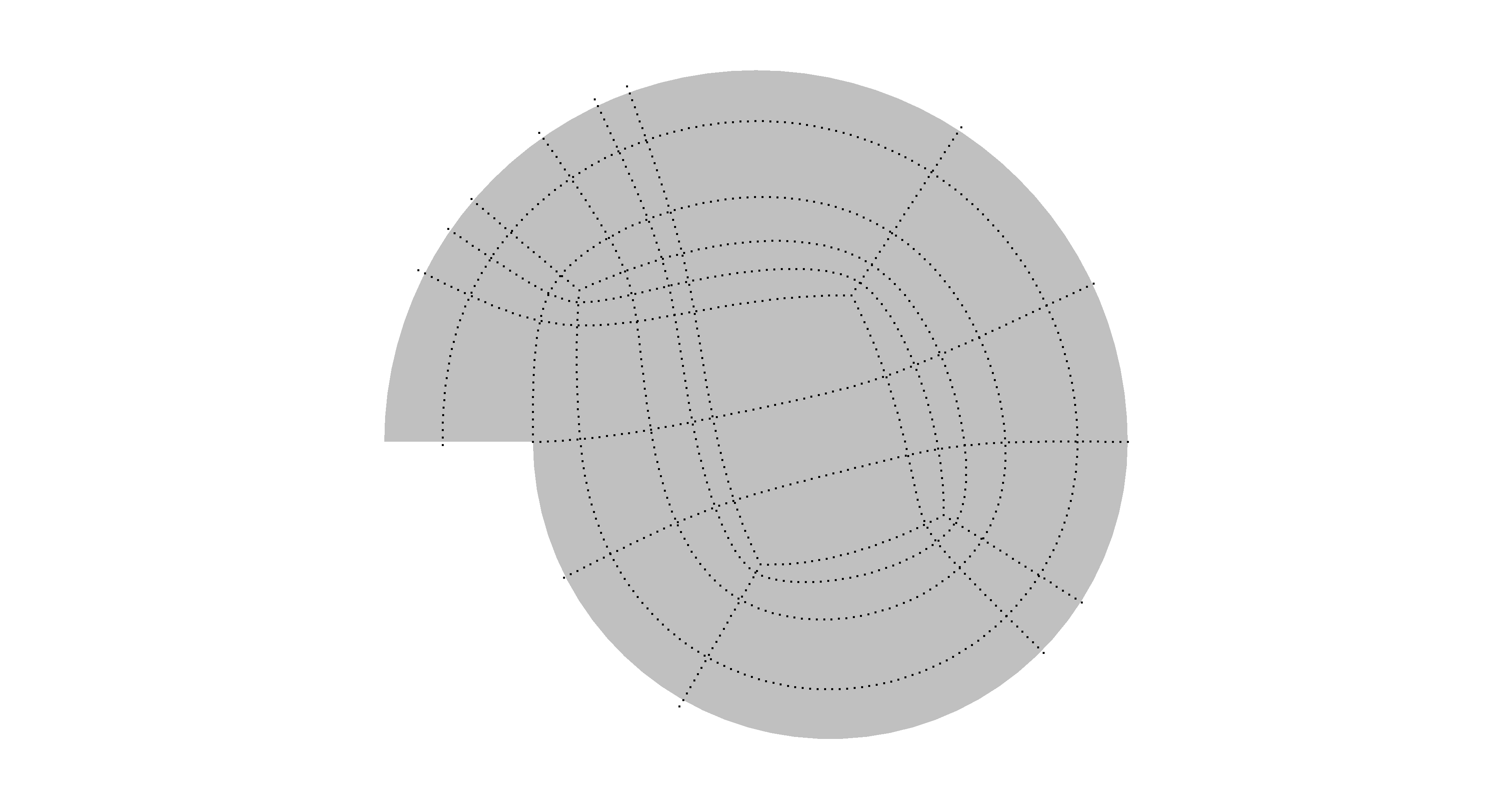}
  \includegraphics[trim={30cm 0 30cm 0},clip,width=\linewidth]{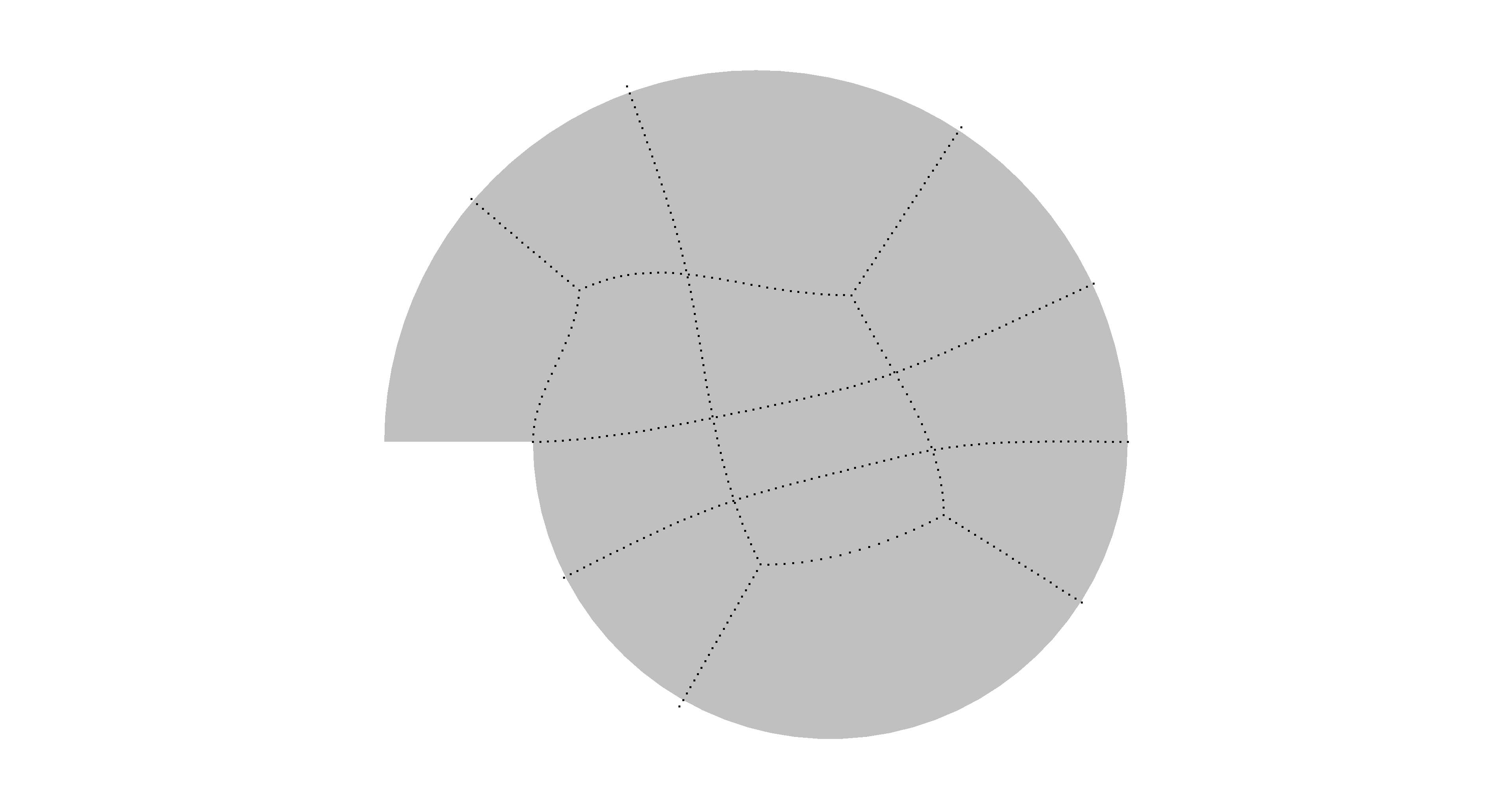}
  \caption{Example of normal (top) and aggressive (bottom) merging. A distance threshold of 5 times the step size was used for aggressive merging.}\label{fig:spiral-streamlines}
\end{figure}

Subsequently, all streamlines are \emph{synchronously} advanced throughout the domain by integrating~\eqref{eq:streamline-integration}.
At any point, two streamlines may meet head on.
For that purpose, we check the proximity of all streamline front points as well as their direction.
If two front points are within merging distance and if they are pointing in opposite directions, their respective streamlines are merged.
This merging distance is typically taken as the integration step size but it can also be artificially increased to perform \emph{aggressive} merging, therefore reducing the total number of streamlines and simplifying the separatrix graph.
Fig.~\ref{fig:spiral-streamlines} illustrates the difference between normal (top) and aggressive (bottom) merging.
When two streamlines are flagged for merging, they are advanced up to each other's starting point and merged pointwise using trigonometric weight functions to retain an optimal streamline direction.
These trigonometric weight functions deflect streamlines as is illustrated in~Fig.~\ref{fig:spiral-streamlines}~(bottom).
Reference~\cite{Marcon} contains more technical details regarding the procedure.

A \(4^{th}\) order multi step Adams-Bashforth integrator is used in physical space.
This integrator requires a search for the element that contains every new point, an operation that is expensive on large meshes.
An inverse mapping is used to transform to parametric coordinates before the high order interpolation fo \(\vec v\) is computed via~\eqref{eq:SolutionInterpolant}.
This process can be accelerated by integrating streamlines in reference space in a similar fashion to the work of Coppola~\textit{et~al.}~\cite{Coppola2001}.
A standard Runge-Kutta scheme can then used while retaining the speed of reference space based integration.

\subsubsection{Midpoint Division}

\begin{figure}[htb!]
  \centering
  \includegraphics[width=\linewidth]{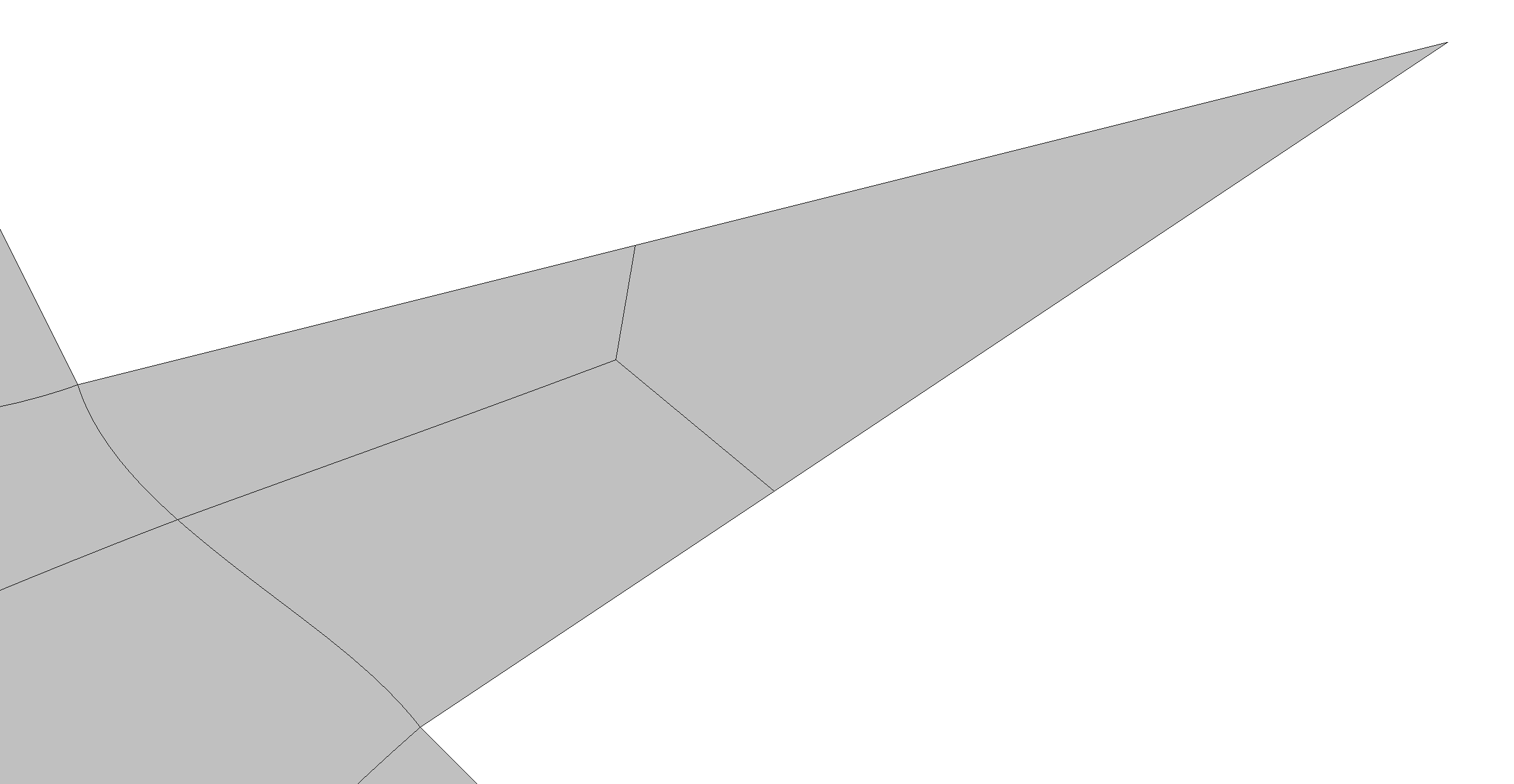}
  \caption{Detail of a domain showing a midpoint division of a degenerate quadrilateral caused by a sharp angled corner. The full domain is shown in Fig.~\ref{fig:polygon-mesh}.}\label{fig:midpoint-division}
\end{figure}

A note should be made about corners with very shallow angles.
Eq.~\eqref{eq:math_corner} can lead to a valence of zero when \( \Delta\theta \) and \( \Delta\psi \) cancel each other.
This represents the situation of all streamlines converging towards a corner.
We call this corner \emph{degenerate} as it will produce a triangular block, i.e.~a collapsed quadrilateral block.
This topology is not valid for the purpose of quadrilateral-only mesh generation, so we resort to \textit{ad hoc} manipulations.
The triangular block can easily be split into three quadrilateral blocks using a midpoint division approach~\cite{Li1995}.
Practically, a 3-valence irregular node is introduced at the center of the triangular block \textit{a posteriori}.
Two of the three streamlines can be extended towards the domain boundaries in a straight line while the third streamline is integrated away from the degenerate corner towards the inside of the domain.
An example of an \textit{ad hoc} midpoint division manipulation is shown in Fig.~\ref{fig:midpoint-division}.

\subsection{Generation of a Quadrilateral Mesh}\label{sec:quad-meshing}

Once a separatrix graph has been generated, it is passed to \emph{NekMesh}, \emph{Nektar++}'s high order mesh generation utility~\cite{Moxey2019}.
\emph{NekMesh} uses OpenCASCADE~\cite{OpenCascadeSAS2018} as its CAD system for high order mesh generation.
The OpenCASCADE kernel offers not only CAD query capabilities but also CAD modification tools.
The discrete separatrices computed by \emph{FieldConvert} are converted into interpolating splines.
These splines can then be used to split the original CAD domains into the desired quadrilateral blocks.
This results in a CAD model of topologically connected quadrilateral faces.
A coarse high order quadrilateral mesh can then be trivially generated by \emph{NekMesh}.

\begin{figure}[htb!]
  \centering
  \includegraphics[trim={5cm 0 5cm 0},clip,width=\linewidth]{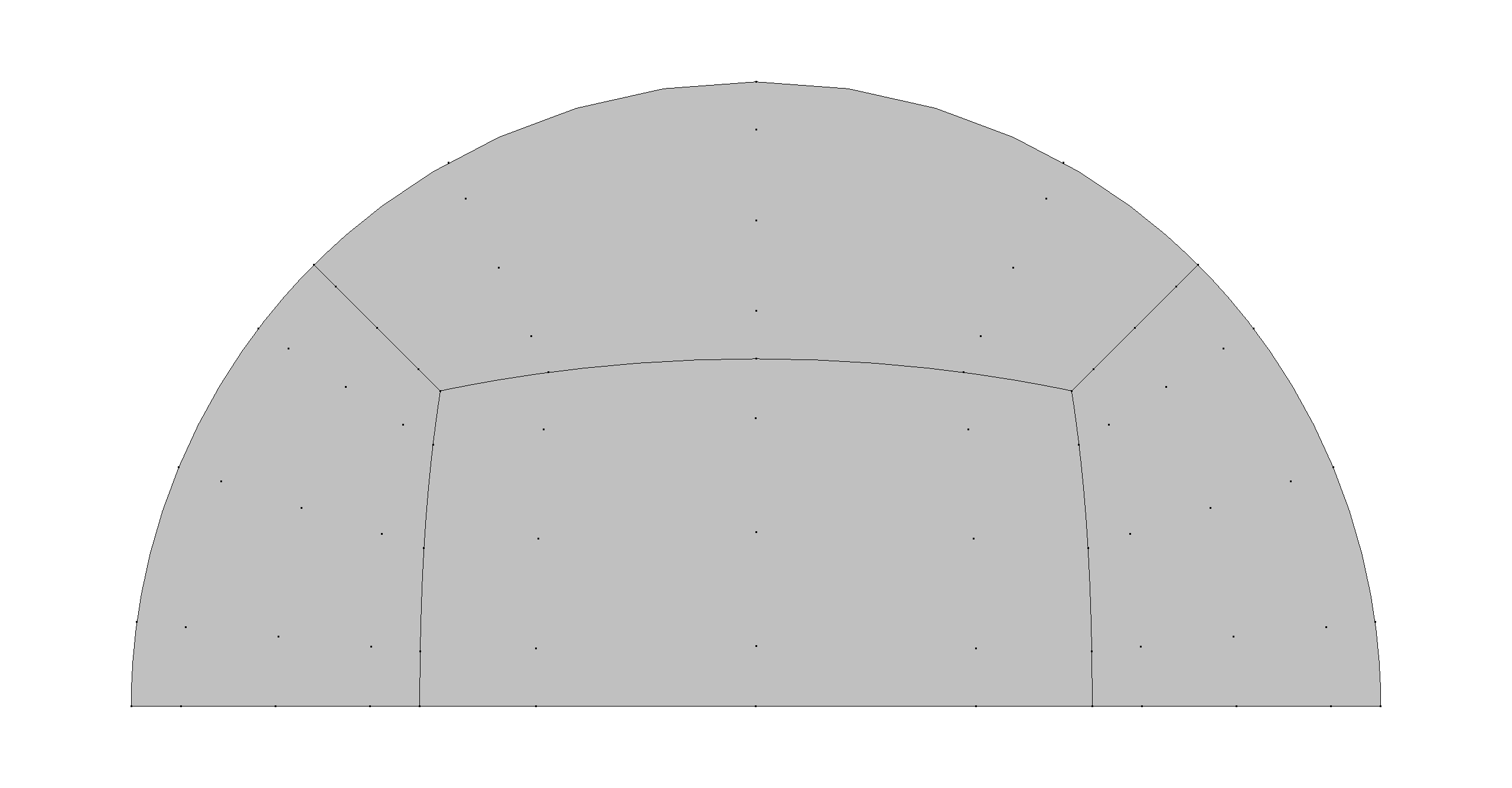}
  \includegraphics[trim={5cm 0 5cm 0},clip,width=\linewidth]{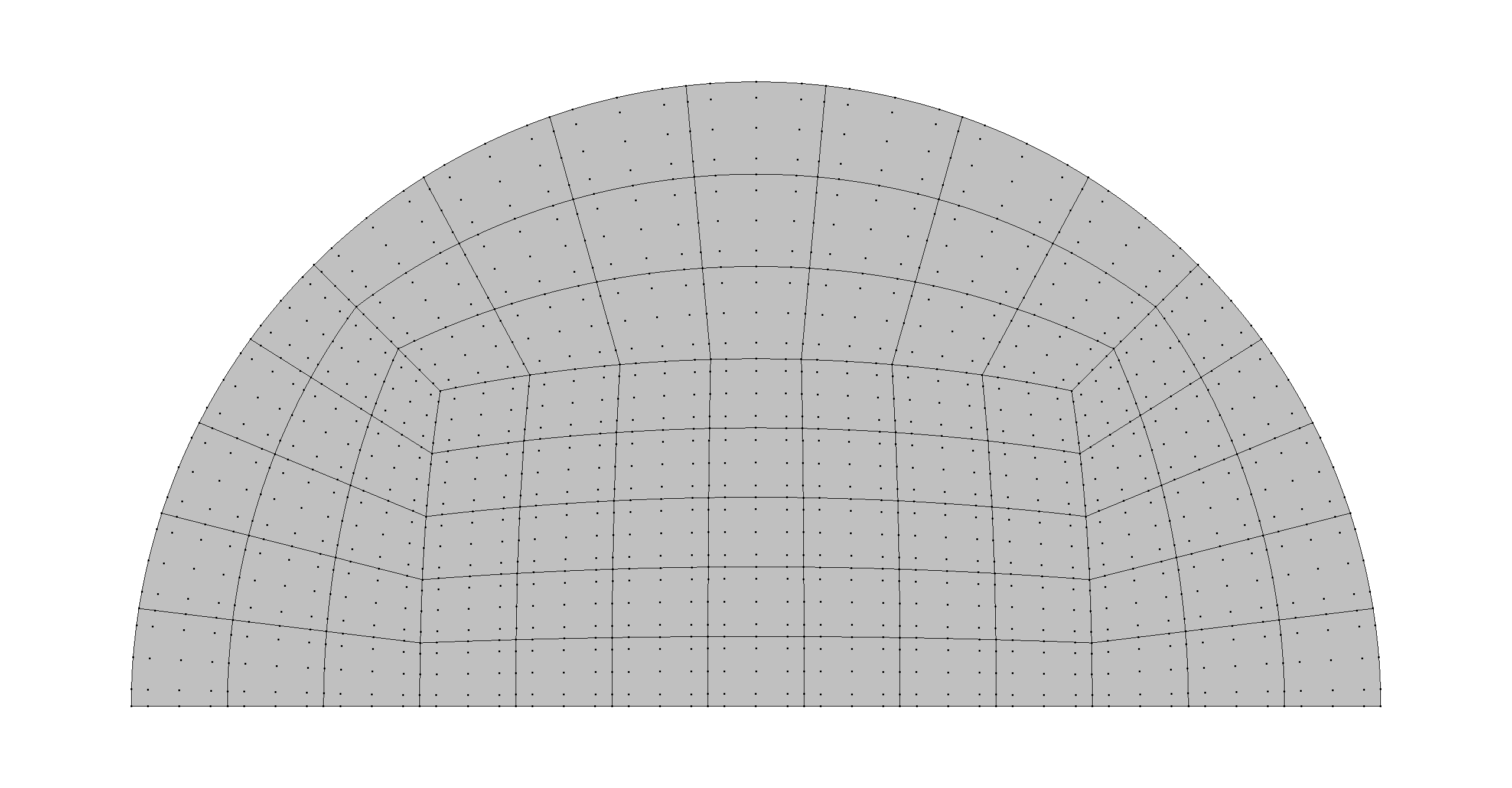}
  \caption{Coarse (top) and subdivided (bottom) quadrilateral meshes obtained on the half disc geometry based on detected critical points and traced and merged streamlines.}\label{fig:half-disc-meshes}
\end{figure}

Such meshes are typically too coarse for spectral element analysis; it is often necessary to split them further.
We take advantage of an isoparametric splitting approach described in~\cite{Moxey2015} to refine strings of blocks, extended in references~\cite{Marcon2018,Marcon2019} for bidirectional splitting.
The method leverages the mapping between physical and reference space coordinates by splitting elements in reference space to obtain valid high order elements in physical space.
In this procedure, it also preserves optimal alignment of blocks.
Fig.~\ref{fig:half-disc-meshes} shows the coarse and fine meshes obtained on the reference geometry.

\section{Advantages of the Method}\label{sec:examples}

We first illustrate advantages of the method with examples before analysing the use of \textit{p}-adaptation in the procedure.
These example geometries, e.g.\ the reference geometry of Fig.~\ref{fig:half-disc-geometry}, are typically used as examples in the cross field literature~\cite{Viertelabs-1708-02316,Nicolas-Kowalski:2012fu,Viertel:2017rt}.
The advantages that we want to illustrate in this section include minimal number of irregular nodes and possibility to aggressively merge streamlines (Geometry \rm{I}), validity of the block decompositions and absence of limit cycles (Geometry \rm{II}), a discretization consistent way of handling corners (Geometry \rm{III}), {\color{black}the preservation of symmetries and patterns (Geometry \rm{IV}),} and finally the flexibility in generating finer meshes (Geometry \rm{V}).

\begin{figure}[htb!]
  \centering
  \includegraphics[trim={8cm 0 10cm 0},clip,width=\linewidth]{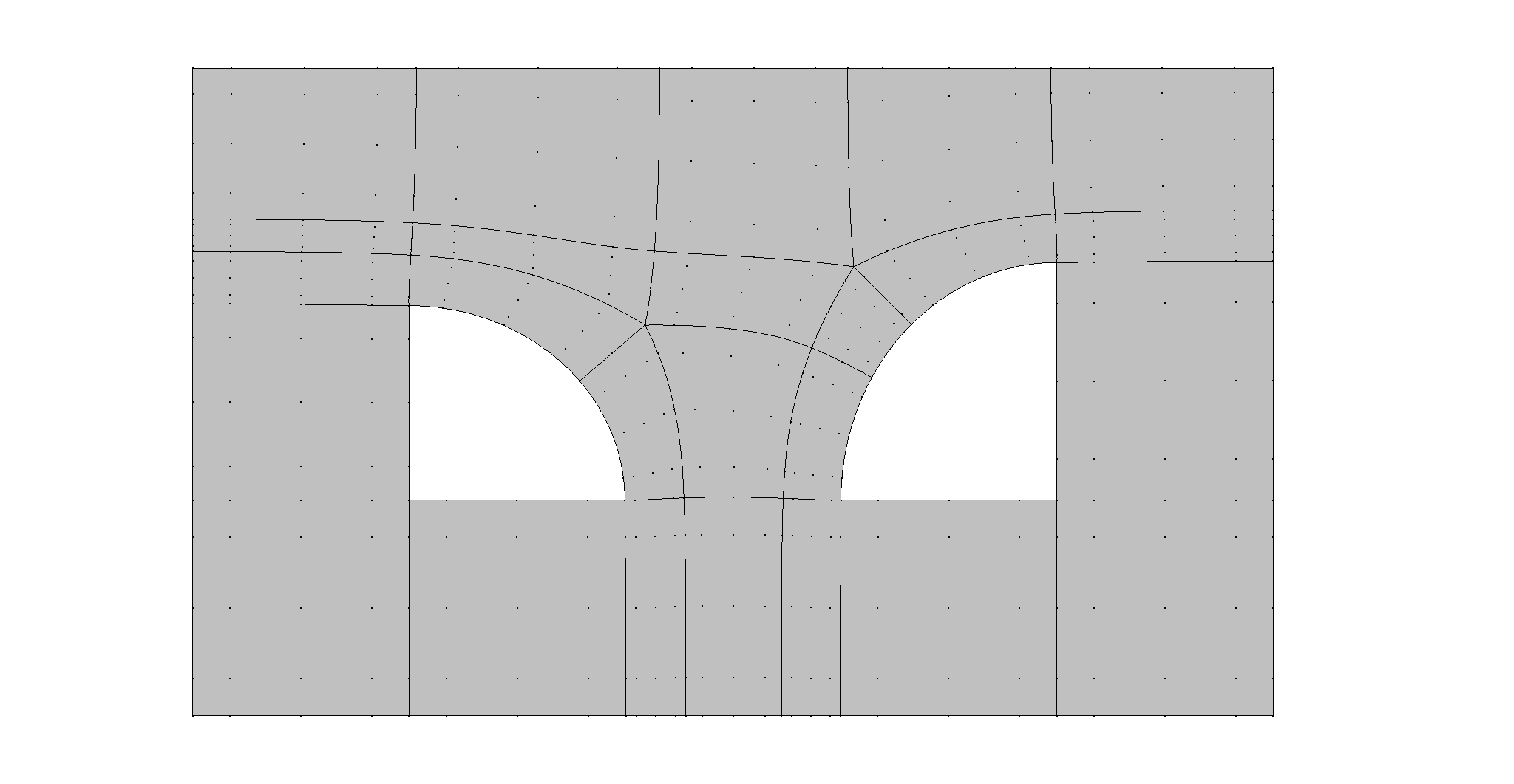}
  \includegraphics[trim={8cm 0 10cm 0},clip,width=\linewidth]{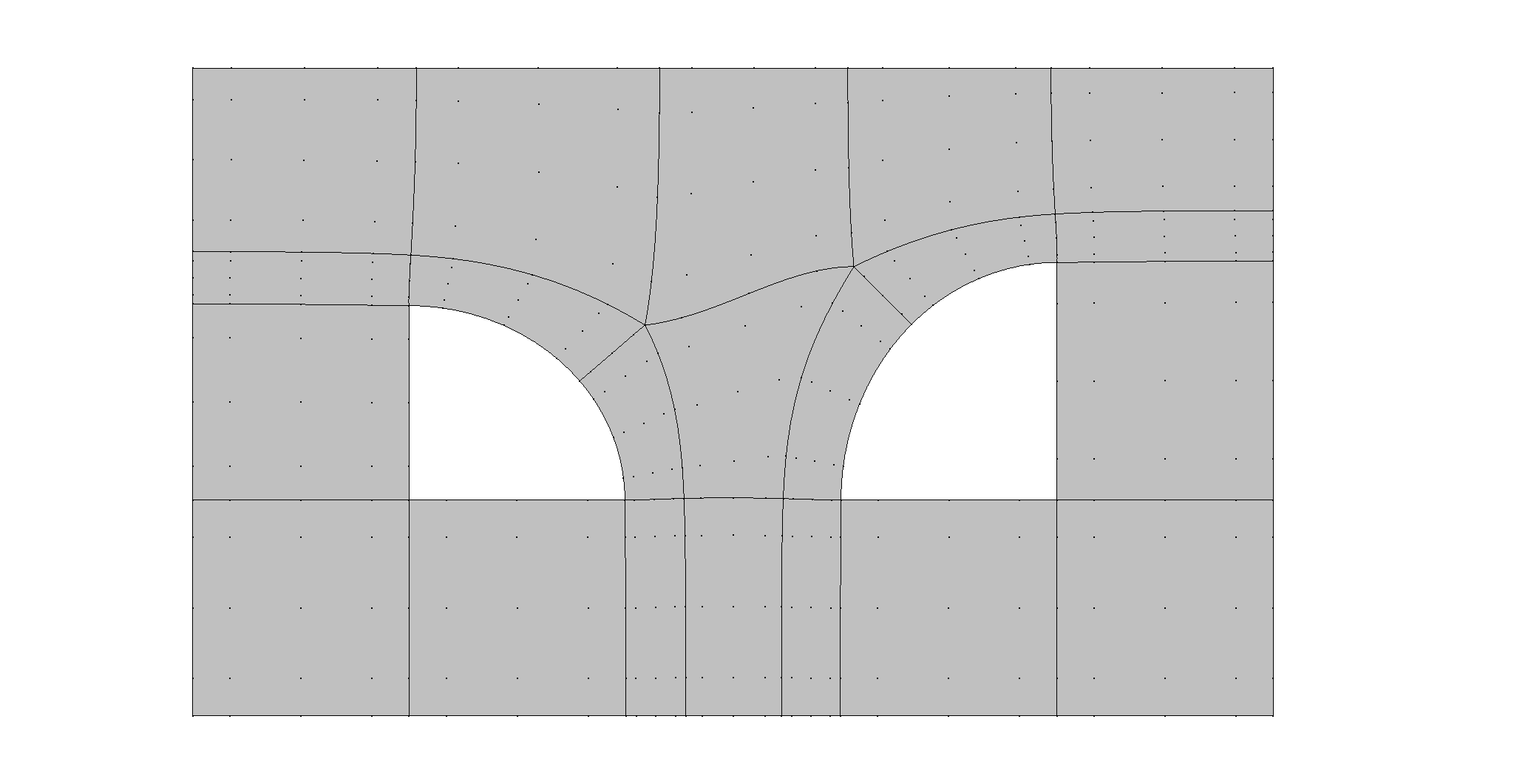}
  \caption{Coarse quadrilateral mesh obtained on Geometry \rm{I} with normal (top) and aggressive (bottom) merging.}\label{fig:iti-mesh}
\end{figure}

Geometry \rm{I} (Fig.~\ref{fig:iti-mesh}) is a multiply connected domain previously meshed using cross fields in~\cite{Nicolas-Kowalski:2012fu}.
While some methods generate spurious irregular nodes in the upper central part of the domain, the approach described here generates only two 3-valence nodes points to obtain an optimal block decomposition of the domain.
Aggressive merging can also be applied to connect the two irregular nodes.
The reader will observe that optimal angles are preserved at the irregular nodes thanks to the trigonometric weights for merging.

\begin{figure}[htb!]
  \centering
  \includegraphics[trim={20cm 0 20cm 0},clip,width=\linewidth]{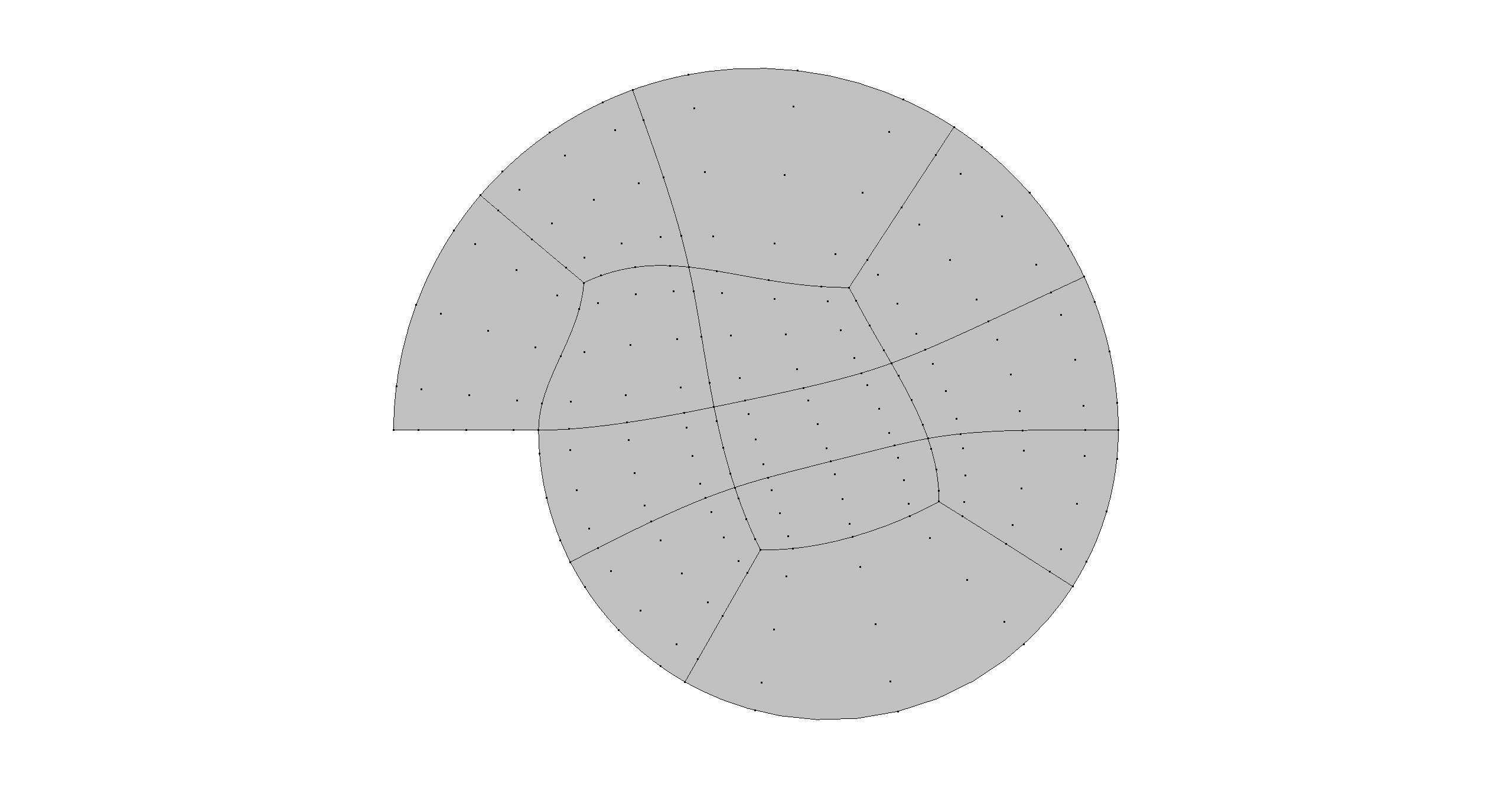}
  \caption{Coarse quadrilateral mesh obtained on Geometry \rm{II}, a nautilus.}\label{fig:spiral-mesh}
\end{figure}

Geometry \rm{II} (Fig.~\ref{fig:spiral-mesh}), sometimes called the \emph{nautilus}, is typically prone to a limit cycle of the streamlines in low order methods leading to invalid block decompositions without \textit{ad hoc} manipulations~\cite{Viertelabs-1708-02316,Viertel:2017rt}.
The formulation of the problem generates a guiding field in which no limit cycle appears.
In fact, the irregular node placement is similar to that of a simple disc, with four 3-valence irregular nodes the optimal number for a valid block decomposition.
As noted in Sec.~\ref{sec:streamlines}, this improved separatrix graph was subject to aggressive streamline merging, for a visually more pleasing coarse mesh.

\begin{figure}[htb!]
  \centering
  \includegraphics[trim={10cm 0 10cm 0},clip,width=\linewidth]{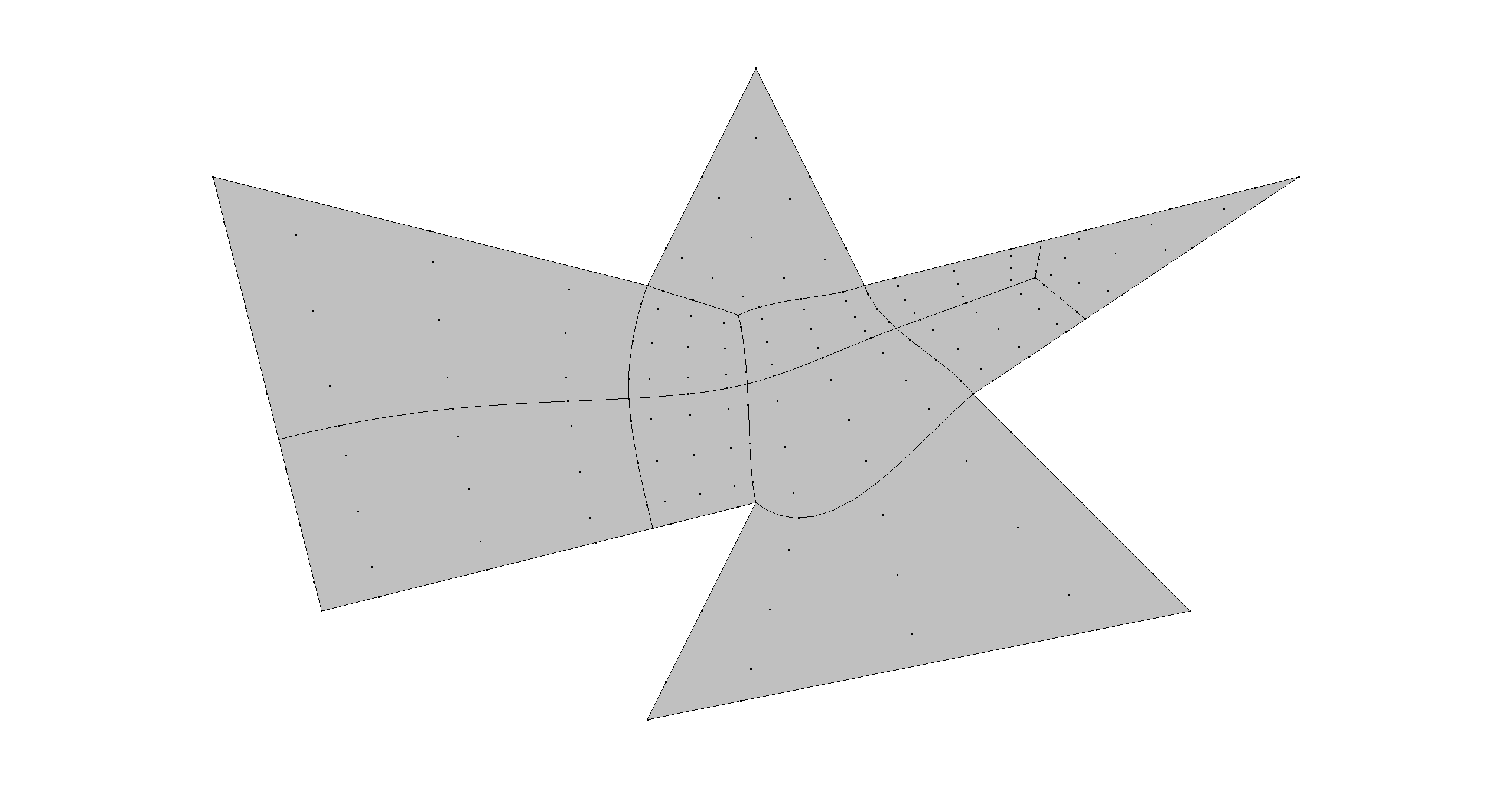}
  \caption{Coarse quadrilateral mesh obtained on Geometry \rm{III}, a polygon.}\label{fig:polygon-mesh}
\end{figure}

Geometry \rm{III} (Fig.~\ref{fig:polygon-mesh}) illustrates the use of a DG discretization to handle angles that are not multiples of \(\pi/2\) without the need for \textit{ad hoc} smoothing of corners.
The rest of the procedure remains however unchanged for easier code implementation and maintenance.
This example also features a sharp corner with a valence of zero whose connected degenerate block is split using the midpoint division method.

{\color{black}
Geometry \rm{IV} (Fig.~\ref{fig:gear-mesh}) illustrates the preservation of symmetries and patterns.
The geometry is that of a simple gear and we observe the same irregular node pattern as seen in the cross field literature~\cite{Viertelabs-1708-02316}.

\begin{figure}[htb!]
  \centering
  \includegraphics[trim={2cm 0 2cm 0},clip,width=\linewidth]{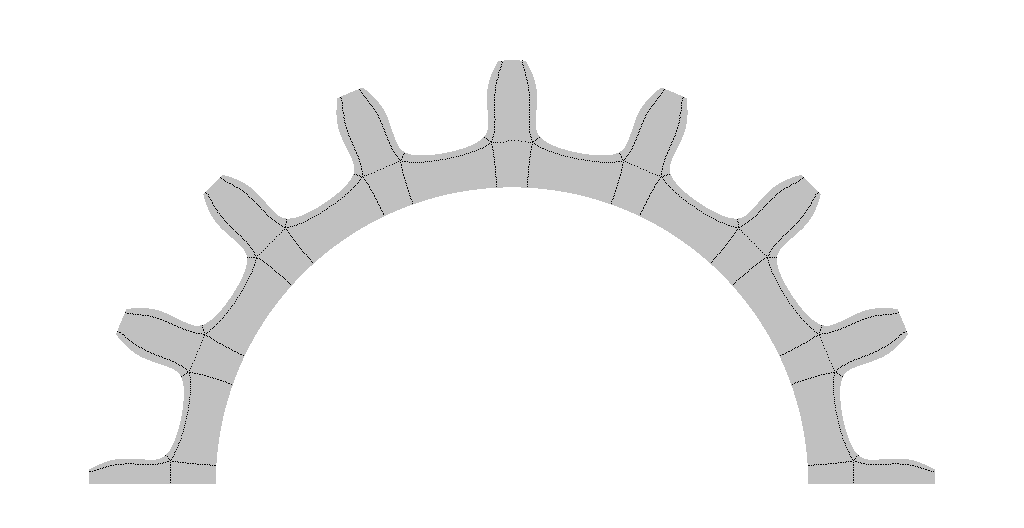}
  \caption{{\color{black}Block decomposition obtained on Geometry \rm{IV}, a gear.}}\label{fig:gear-mesh}
\end{figure}
}

\begin{figure}[htb!]
  \centering
  \includegraphics[trim={15cm 0 15cm 0},clip,width=\linewidth]{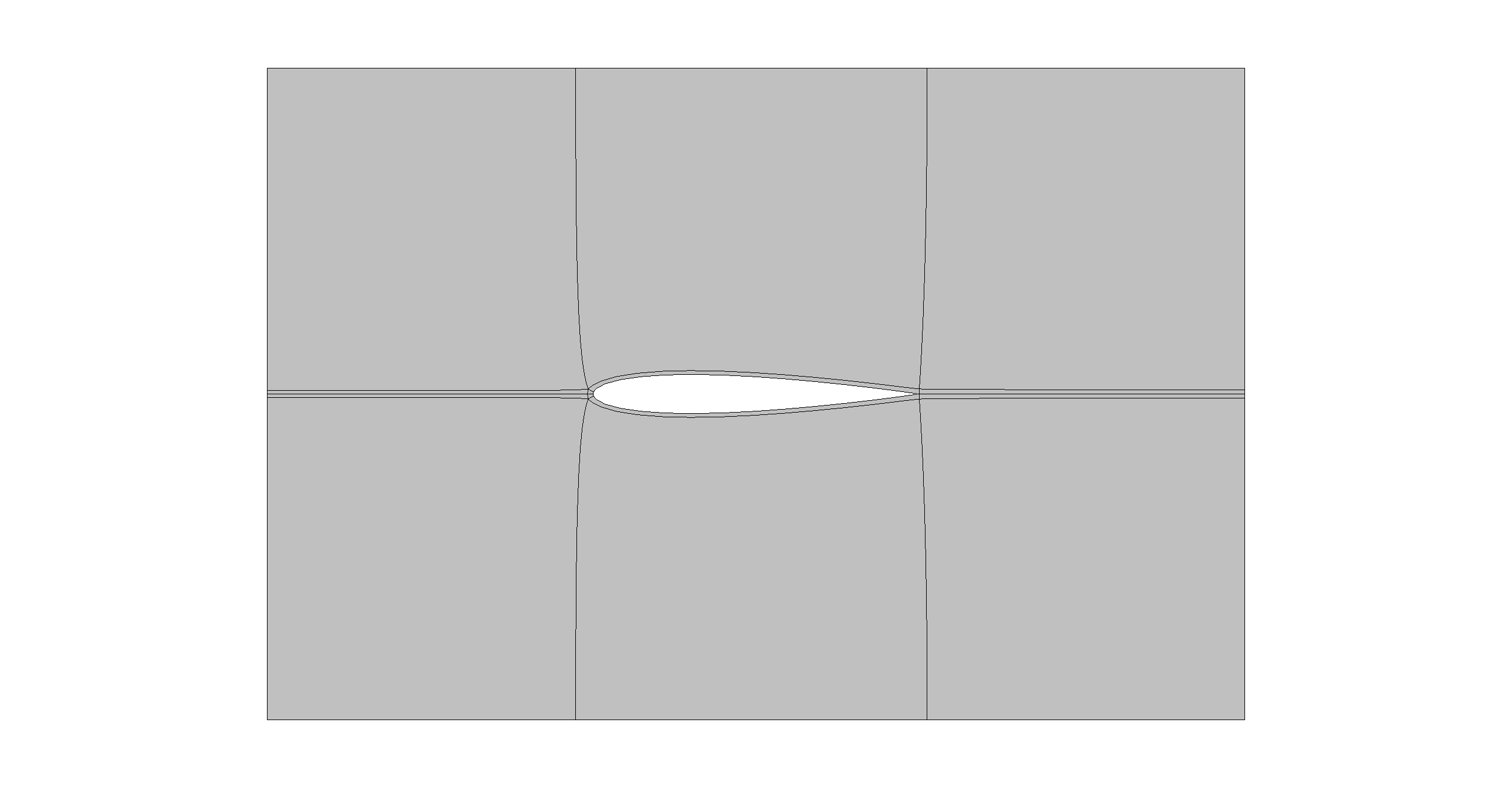}
  \includegraphics[trim={15cm 0 15cm 0},clip,width=\linewidth]{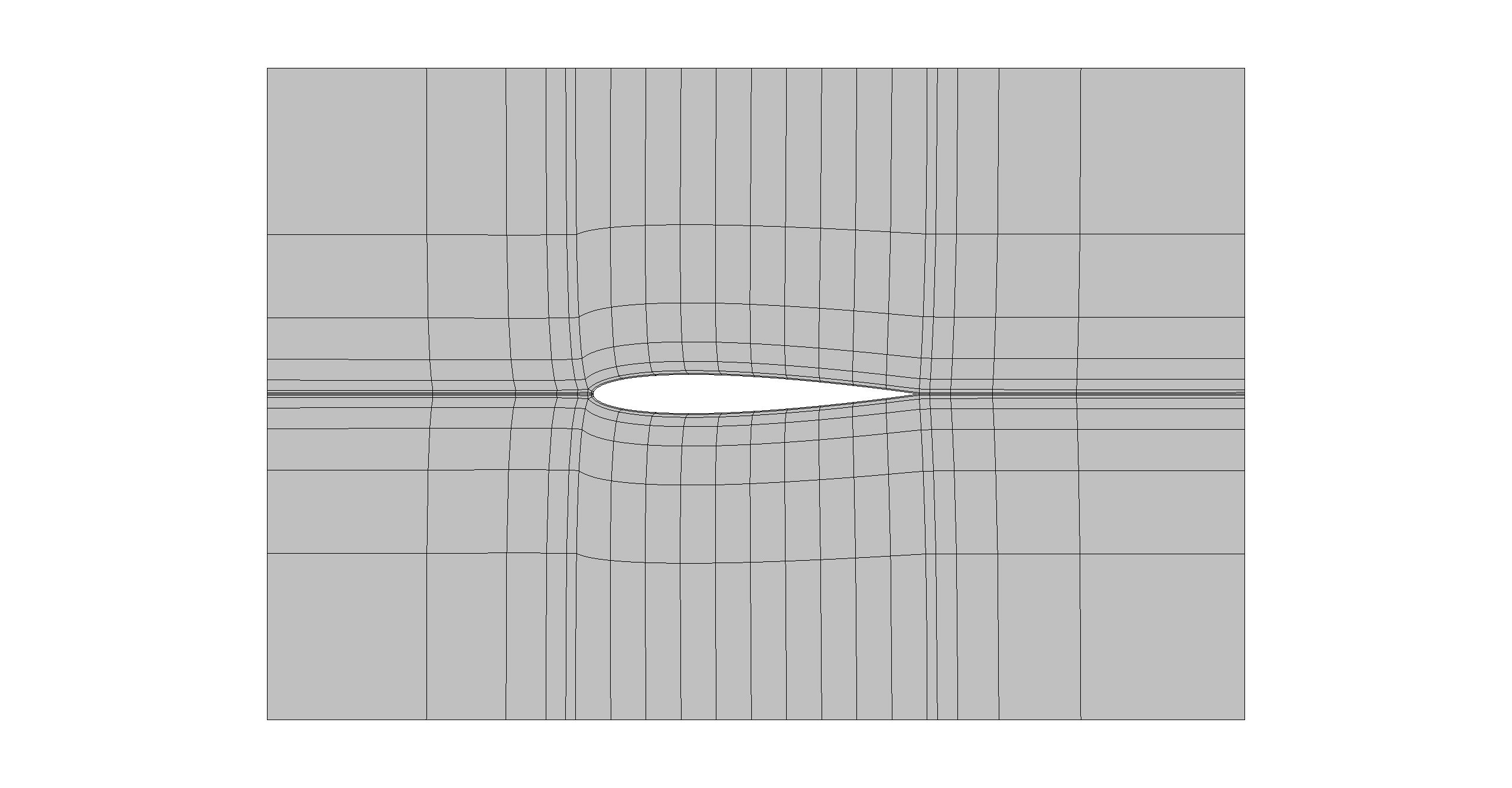}
  \caption{Coarse (top) and split (bottom) quadrilateral meshes obtained on Geometry \rm{V}, a NACA 0012 profile. High order interior degrees of freedom have been hidden for clearer visualization.}\label{fig:naca-mesh}
\end{figure}

Whereas Geometries \rm{I}--\rm{III} illustrate the generation of naturally curved very coarse meshes, finer meshes are typically required for SEM analysis with user defined element size requirements.
Geometry \rm{V} (Fig.~\ref{fig:naca-mesh}) illustrates the use of the isoparametric splitting to obtain a finer mesh with a certain distribution of elements.
The geometry is that of NACA 0012 aerofoil, a typical geometry used in academic CFD studies.

\begin{figure}[htb!]
  \centering
  \includegraphics[trim={15cm 0 15cm 10cm},clip,width=\linewidth]{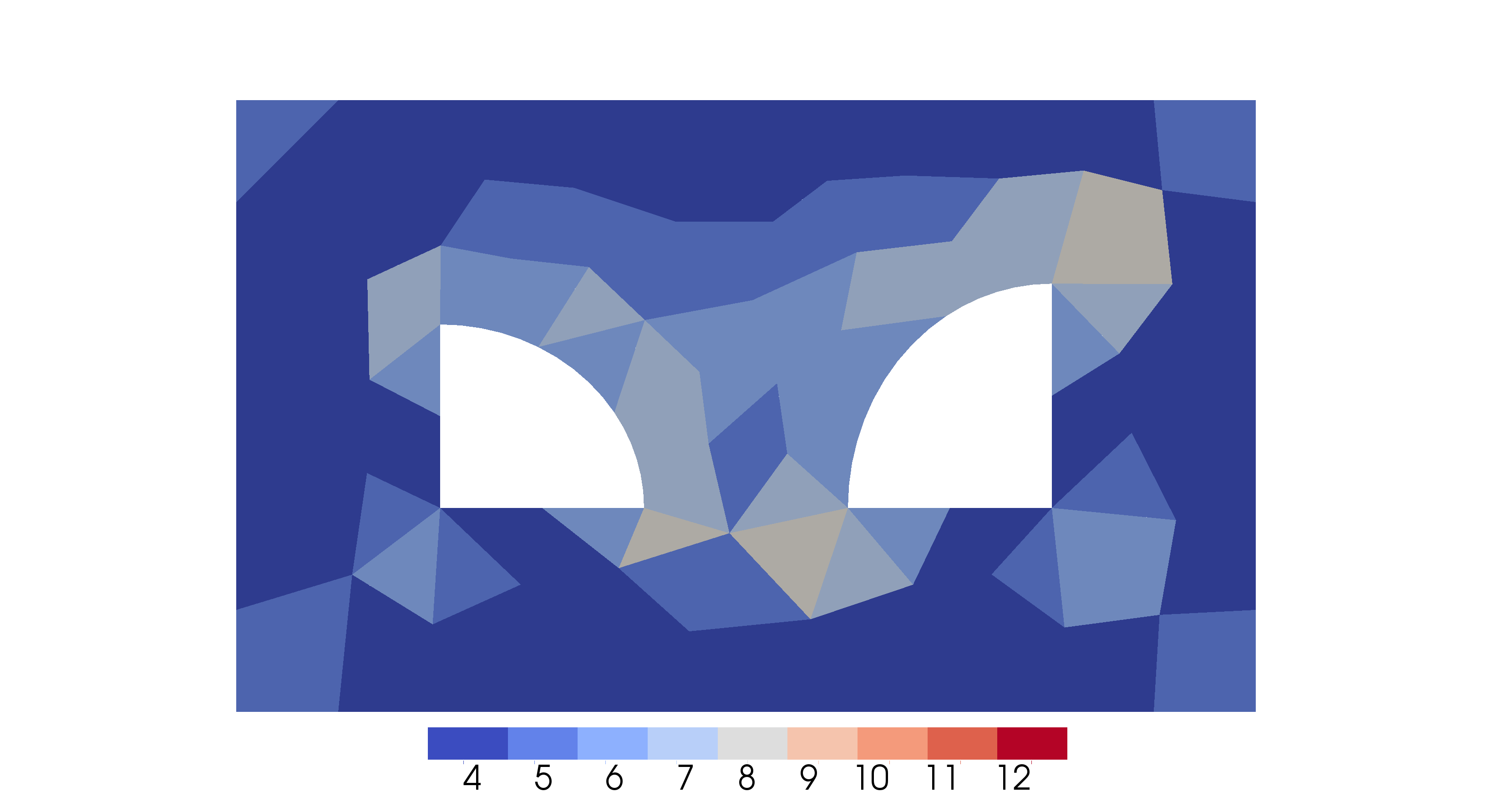}
  \includegraphics[trim={15cm 0 15cm 10cm},clip,width=\linewidth]{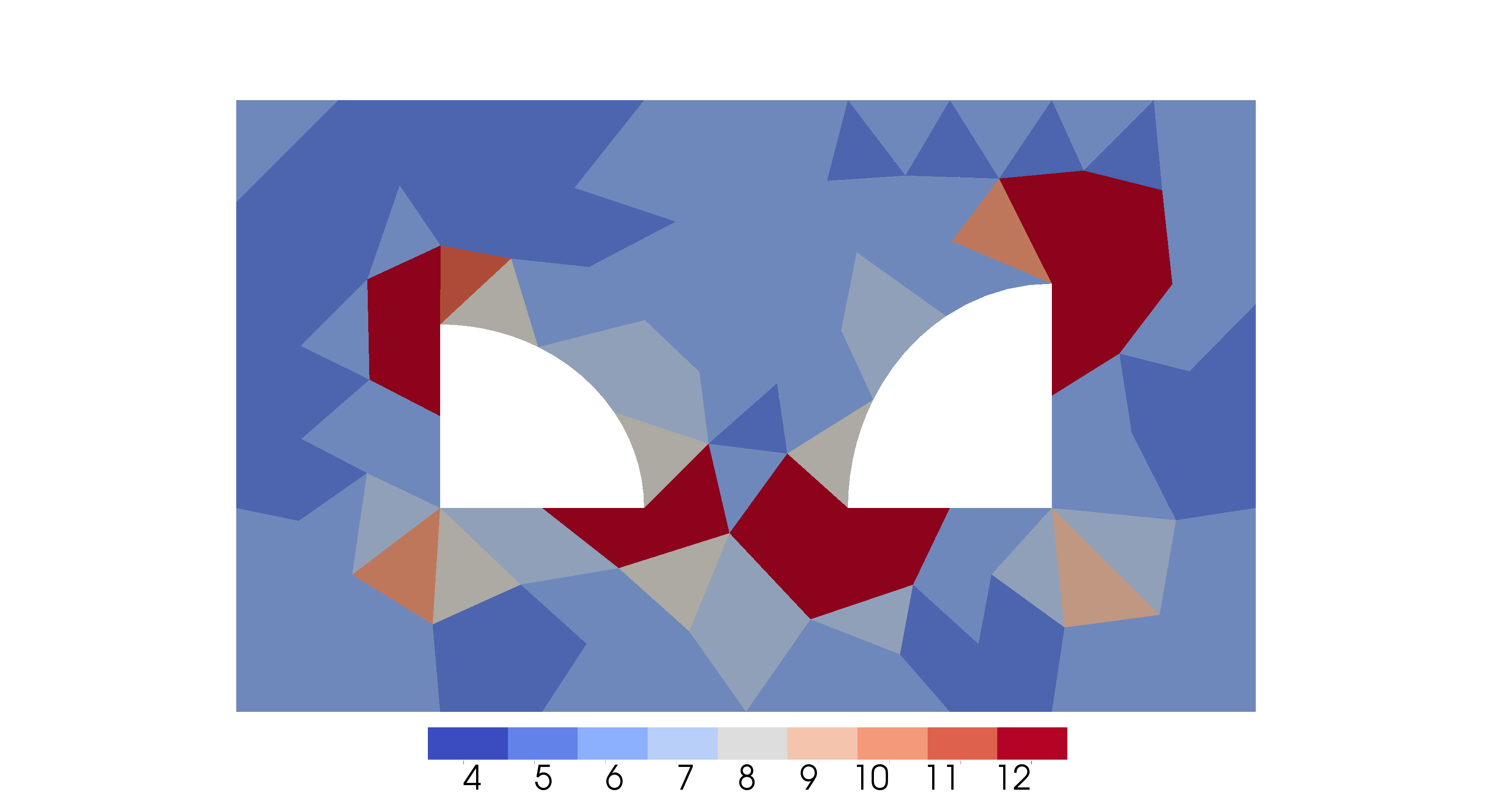}
  \caption{Distribution of local number of modes (\(=P+1\)) for Geometry \rm{I}. Adaptation based on \(u\) (top) and \(v\) (bottom).}\label{fig:iti-p}
\end{figure}

We next demonstrate the use of \textit{p}-adaptation.
For Geometry \rm{I}, we show the distribution of the local number of modes (\(=P+1\)) after \textit{p}-adaptation in Fig.~\ref{fig:iti-p}.
The local polynomial order is higher mostly near boundaries and in the central area where rapid change in the solution is observed.
The local polynomial order also seems unexpectedly high in the \(v\)-based case but this can be easily explained by the actual solution.
Because the geometry is dominated by {\color{black}Cartesian}-aligned straight-sided boundaries, a large part of the field (especially the areas in the bottom left and right) has values close to {\color{black}\( \vec v =  (1,0) \)}.
This means that the low order modes have no energy {\color{black}for the \(v\) component}.
{\color{black}Therefore, for the \(v\) based sensor,} any energy appearing in the higher order modes will trigger an increase in the local polynomial order.
For such examples dominated by {\color{black}Cartesian}-aligned boundaries, it is therefore best to base the sensor on \(u\) {\color{black}since \(u_b = 1\) in Cartesian directions}.
The same phenomenon is retrospectively observed for the half disc geometry in Fig.~\ref{fig:half-disc-p}.

{\color{black}
The use of high order methods might intuitively sound prohibitively expensive for this method.
We find, however, that the computational time remains low, with all examples in this paper taking \(\mathcal{O}(10^{0})\) seconds to process.
The main cost in \textit{NekMesh} is the computation of the streamlines because that is currently done in physical, not reference, space.
Streamline integration currently takes between roughly 50\% and 90\% of the time of the whole procedure, depending on the size of the domain and of the base mesh, the integration step size, and how far streamlines have to travel.
All other tasks --- each solution of the Laplace problem, the analysis of the critical points and corners and the generation of a quad mesh --- take about the same share of the remaining 10\% to 50\% of computational time.
}

\section{Summary and Conclusions}

We have demonstrated the benefits of a field guided method to generate quadrilateral block decompositions, and meshes, of multiply connected two dimensional domains.
The guiding field is inspired by cross fields but crosses are, in fact, never generated.
This method generates naturally curved elements suitable for spectral element methods, and other related methods, and it does so by using a spectral element solver.
The guiding field is obtained from a boundary value problem using a Laplace solver, with boundary conditions determined by local geometrical information.
This information can be discontinuous, in which case a discontinuous Galerkin formulation is used for a discretization consistent way of imposing these discontinuous boundary conditions.
The separatrix graph, obtained from the analysis of the guiding field, is later used to split the two dimensional model.
From this block decomposed geometry, a coarse, valid, curved mesh is generated that can be refined to users' needs.
The procedure is implemented in the open source spectral/\textit{hp} element framework \textit{Nektar++} and its high order mesh generation tool \textit{NekMesh}.

We emphasize some of the benefits of this method including some advantages over traditional cross field based approaches as well.
Whereas high order mesh generation techniques typical suffer from invalid elements, the method presented here creates coarse \textit{a priori} curved meshes that are naturally valid.
It does so by using a spectral element solver, which not only benefits from spectral convergence but is also suitable to \textit{p}-adaptation for better computational resources.
This creates a high resolution guiding field, due to the high order information of the solution, used for the creation of an optimal separatrix graph: small number of irregular nodes, preservation of symmetries, absence of limit cycles, etc.
The analysis of the guiding field is performed in reference space at a low computational cost.
The refinement of the coarse mesh is later performed using an isoparametric splitting approach that preserves the high order validity of the mesh.

The authors envision two improvements of this proof of concept for better relevance in the CFD community.
On the one hand, an extension to three dimensional hexahedral block decompositions must be sought.
This extension requires a reformulation of the problem based on three dimensional frame arithmetic.
Promising work of Ray~\textit{et~al.}~\cite{Ray2016} is paving the way in that direction.
On the other hand, the generation of metric aligned quadrilateral blocks, as was done in reference~\cite{Fogg:2013kl}, seems possible through the use of a forcing term.
The solution of a Poisson problem with adequate forcing would enable modifications of the separatrix graph with respect to a metric field.
Such a solution would also be appropriate for the relocation of irregular nodes away from boundaries with strong curvature (see for example Fig.~\ref{fig:naca-mesh}).

\section*{Acknowledgements}

This project has received funding from the European Union's Horizon 2020 research and innovation programme under the Marie Sk\l{}odowska-Curie grant agreement No 675008.
This work was supported by a grant from the Simons Foundation (\#426393, David Kopriva).
The first two authors would also like to thank Prof.~Gustaaf Jacobs of the San Diego State University for his hospitality.
Finally, JM would like to thank the committee of the 28\(^{th}\) International Meshing Roundtable for their travel financial support.

\bibliography{refs}

\end{document}